\documentclass{amsart}%
\usepackage{amsfonts}
\usepackage{amsmath}
\usepackage{amssymb}
\usepackage{graphicx}
\usepackage{hyperref}%
\providecommand{\U}[1]{\protect\rule{.1in}{.1in}}
\newtheorem{theorem}{Theorem}
\theoremstyle{plain}

\newtheorem{case}{Case}
\newtheorem{claim}{Claim}

\newtheorem{corollary}[theorem]{Corollary}

\newtheorem{definition}{Definition}
\newtheorem{example}[theorem]{Example}

\newtheorem{lemma}[theorem]{Lemma}

\newtheorem{remark}[theorem]{Remark}

\numberwithin{equation}{section}
\numberwithin{theorem}{section}
\makeatletter
\@addtoreset{case}{theorem}
\@addtoreset{claim}{theorem}
\makeatother
\begin{document}
\title[Invariant subspaces]{Invariant subspaces of subgraded Lie algebras of compact operators}
\author{Matthew Kennedy}
\address{Department of Pure Mathematics, University of Waterloo, 200 University Avenue
West \\
Waterloo, Ontario, Canada N2L 3G1 }
\email{m3kennedy@uwaterloo.ca}
\author{Victor S. Shulman}
\address{Department of Mathematics, Vologda State Technical University \\
15 Lenina Street, Vologda 16000, Russian Federation}
\email{shulman\_v@yahoo.com}
\author{Yuri V. Turovskii}
\address{Institute of Mathematics and Mechanics, National Academy of Sciences of Azerbaijan\\
9 F. Agayev Street, Baku AZ1141, Azerbaijan}
\email{yuri.turovskii@gmail.com}
\thanks{}
\thanks{}
\thanks{The support received from INTAS project No 06-1000017-8609 is gratefully
acknowledged by the third author}
\thanks{This paper is in final form and no version of it will be submitted for
publication elsewhere.}
\date{}
\subjclass[2000]{Primary 47L70, 47A15; Secondary 47L10, 16W25}
\keywords{Graded Lie Algebra, Engel Lie Algebra, Lie Triple System, Jordan Algebra,
Compact Banach Algebra, Compact Operator, Invariant Subspace,
Triangularization, Cartan Criterion. }
\dedicatory{ }
\begin{abstract}
We show that finitely subgraded Lie algebras of compact operators have
invariant subspaces when conditions of quasinilpotence are imposed on certain
components of the subgrading. This allows us to obtain some useful information
about the structure of such algebras. As an application, we prove a number of
results on the existence of invariant subspaces for algebraic structures of
compact operators, in particular for Jordan algebras and Lie triple systems of
Volterra operators. Along the way we obtain new criteria for the
triangularizability of a Lie algebra of compact operators.

\end{abstract}
\maketitle

\section{Introduction}

\subsection{Aim of the work and remarks on the history of the topic}

The classical Invariant Subspace Problem, the problem of the existence of a
nontrivial closed subspace invariant under an operator from a given class,
extends in various forms to families of operators. Even partial positive
results on this subject can be of considerable importance for the analysis of
various algebraic systems (groups, semigroups, and both associative and
non-associative algebras) of operators, and for the structure and
representation theory of abstract topological-algebraic systems. Typically,
the application of these results involves the construction of a maximal chain
of invariant subspaces for a system of operators. This process is often
referred to as \textit{triangularization}, and in fact, it provides many of
the same advantages that triangular forms give for the analysis of systems of matrices.

There are two popular versions of the invariant subspace problem. The first
concerns commutative families of operators, and the second general operator
algebras (i.e. problems of Burnside type). The solution of both versions for
the class of compact operators was given in 1973 by Victor Lomonosov
\cite{Lom}, who proved that every algebra of compact operators on a Banach
space $\mathcal{X}$ possessing no invariant subspaces (i.e., which is
\textit{irreducible}) is dense in the algebra $\mathcal{B}(\mathcal{X})$ of
all bounded operators on $\mathcal{X}$, with respect to the weak operator
topology. (Notice this implies that any commutative family of compact
operators has an invariant subspace). An important consequence is that an
algebra of compact quasinilpotent (\textit{Volterra}, for brevity) operators
has invariant subspace (i.e., is \textit{reducible}). This is the core of the
results of \cite{Lom}; the general case can be easily deduced from it by using
spectral projections.

For the investigation and classification of general algebraic systems of
compact operators, the case when the system consists of Volterra operators is
a natural first step. It can be said that such systems play the central role
in the so called \textquotedblleft radical\textquotedblright\ part of the more
general structural theory. Indeed, such an approach was very successful in the
theory of finite dimensional Lie algebras. We refer, for instance, to the
classical theorems of Engel and Lie.

The first attempts to extend the results of Lomonosov to Lie algebras of
compact operators \cite{Woj}, and to semigroups of compact operators
\cite{Rad}, used techniques based on the trace, so positive results were
obtained under the condition that the Lie algebra or semigroup in question
contain a nonzero trace class operator (or more generally, an operator from a
Schatten class). With the introduction of more general techniques based on the
joint spectral radius \cite{Sh83, T84, T85}, more general results, free from
the above restriction were obtained. It was proved in \cite{Tur1998} that a
semigroup of Volterra operators has an invariant subspace, and the
reducibility of Lie algebras of Volterra operators was proved in
\cite{ShT2000}. The latter result allowed for the development of a detailed
theory of the structure of Lie algebras of compact operators in \cite{ShT2005}.

More recently, in \cite{Ken}, the existence of invariant subspaces was
established for Jordan algebras of Volterra operators which contain nonzero
trace class operators. The initial aim of the present work was to extend this
result by avoiding the hypothesis of the existence of nonzero trace class
operators, that is, to prove that every Jordan algebra of Volterra operators
has an invariant subspace. Upon further examination, we realized that the
needed result could be obtained within the much more general framework of the
theory of subgraded Lie algebras of compact operators.

\begin{definition}
\label{grad} Let $\Gamma$ be an abelian group. A Lie algebra $L$ is said to be
$\Gamma$-subgraded if it can be written as a sum of subspaces indexed by the
elements of $\Gamma$,
\[
L=\sum_{\gamma\in\Gamma}L_{\gamma},
\]
which satisfy the condition
\[
\lbrack L_{\gamma},L_{\delta}]\subseteq L_{\gamma+\delta}.
\]
If $L=\oplus_{\gamma\in\Gamma}L_{\gamma}$, i.e. if the sum is direct, then $L$
is said to be $\Gamma$-graded. Elements of $\cup_{\gamma\in\Gamma}L$ are
called homogeneous, and elements of $L_{\gamma}$ are called $\gamma$-homogeneous.
\end{definition}

Unless otherwise specified, the index group $\Gamma$ is assumed to be finite,
and in this case we say that $L$ is \textit{finitely subgraded}.

There are two types of quasinilpotence conditions that we will impose on an
element $a$ of a Lie algebra $L$ of compact operators. The first condition,
which is purely operator-theoretic, is that $a$ is a Volterra (i.e.
quasinilpotent) operator. The second condition, which is less restrictive, is
that the operator
\[
\operatorname*{ad}a:x\rightarrow\lbrack a,x],
\]
defined on $L$, is quasinilpotent. Notice that the second condition makes
sense in the context of any abstract Banach Lie algebra, since
$\operatorname*{ad}a$ is the image of $a$ in the adjoint representation of
$L$. When $a$ satisfies the second condition, we say that $a$ is an
\textit{Engel} element of $L$.

In dealing with subgraded Lie algebras of compact operators, the task we set
for ourselves is to prove their reducibility or triangularizability by
imposing a minimal number of quasinilpotence conditions. In particular, we
will consider two very natural conditions -- that all their homogeneous
elements are Engel, or that all their $0$-homogeneous elements are Engel.

There is a special motivation for considering the case of $\mathbb{Z}_{2}%
$-subgraded Lie algebras with the condition that the $1$-homogeneous elements
are Engel. For, if $J$ is a Jordan algebra of operators then it can be checked
that $L=[J,J]+J$ is a $\mathbb{Z}_{2}$-subgraded Lie algebra, with
$L_{0}=[J,J]$ and $L_{1}=J$. One should note that the intersection of $L_{0}$
and $L_{1}$ may be nonzero.

During our work, we obtained a number of conditions, interesting in their own
right, which imply the reducibility and triangularizability of (not
necessarily subgraded) Lie algebras of compact operators. These conditions
proved valuable in our investigation of subgraded Lie algebras and Jordan
algebras of compact operators.

\subsection{Structure of paper and outline of results}

This paper is naturally partitioned into three interrelated parts, with each
part building up to ``main results.''

\subsubsection{Triangulation criteria for Lie algebras of compact operators}

In Part \ref{__1criteria} (Sections \ref{_2preliminaries}-\ref{_4further}), we
develop a list of criteria for the reducibility and triangularizability of Lie
algebras of compact operators.

Section \ref{_2preliminaries}, being a background for the rest of the paper,
contains a review of the known criterion for reducibility and
triangularizability. Among the most important is an infinite dimensional
analogue of Cartan's Solvability Lemma (see Lemma \ref{carfin}), which is a
variation of a result in \cite{Ken}.

There is a well-known condition \cite{Mu1} for the triangularizability of an
associative algebra of compact operators. Namely, a closed associative algebra
of compact operators is triangularizable if and only if it is commutative
modulo its Jacobson radical. For Lie algebras, a similar result, formulated in
"radical-like" terms, and without referring to any underlying space, was
obtained in \cite{ShT2005}.

Let us call a Banach Lie algebra $L$ \textit{Engel} if all its elements are
Engel, and $E$\textit{-solvable} if every non-zero quotient of $L$ by a closed
ideal of $L$ has a nonzero closed Engel ideal. (When $L$ is finite
dimensional, this condition is equivalent to the solvability of $L$). The
above-mentioned result of \cite{ShT2005} states that a Lie algebra of compact
operators is triangularizable if and only if it is $E$-solvable. Furthermore,
if $L$ has a non-scalar $E$-solvable ideal, then it is reducible.

The main result of Section \ref{_3subalgebra} is that for a Lie algebra $L$ of
compact operators, the Engel elements of $L$ which belong to an $E$-solvable
subalgebra $W$ of $L$ form an ideal of $W$. Under certain additional
conditions, this allows us to obtain an Engel ideal of $L$ itself, and as a
consequence, to prove the reducibility of $L$. In the proof we make use of
techniques from the theory of compact Banach algebras.

In Section \ref{_4further}, we further develop the list of criteria which
imply the reducibility of a Lie algebra of compact operators. Among the most
effective are those that impose restrictions on finite rank or
$\operatorname*{ad}$-nilpotent elements. For example, when $L$ is closed, we
show that if the commutator $[a,b]$ of arbitrary finite rank operators $a,b\in
L$ is nilpotent, then $L$ is triangularizable. This greatly extends the result
of Katavolos and Radjavi \cite{KaRad} that an associative algebra of compact
operators is triangularizable if all its commutators are quasinilpotent.

In general, the first step in establishing the triangularizability of a Lie
algebra $L$ of compact operators is to show that it is reducible.
Unfortunately however, conditions imposed only on finite rank operators are
not necessarily preserved in passing to the (closure of the) Lie algebras
induced by quotients of invariant subspaces, and as a result, the argument
cannot be repeated. This is because it is possible for extra finite rank
operators to appear.

To circumvent this obstacle, we first prove that the closed associative
algebra $\mathcal{A}(L)$, generated by $L$, is commutative modulo its radical
$\operatorname*{rad}\mathcal{A}(L)$, which establishes the triangularizability
of $L$ by the above-mentioned criterion of \cite{Mu1}. The quotient
$\mathcal{A}(L)/\operatorname*{rad}\mathcal{A}(L)$ is not an algebra of
compact operators, but it is a compact Banach algebra. This is why the results
on compact Banach algebras obtained in Section \ref{_3subalgebra} are
important for our work.

For (not necessarily closed) Lie algebras of compact operators, another
particularly important criterion for triangularizability is the condition that
all  commutators are Engel elements (or Volterra operators). In Part
\ref{__2graded}, we extend this result to subgraded Lie algebras of compact
operators, requiring only that  homogeneous commutators are Engel elements (or
Volterra operators).

\subsubsection{Subgraded Lie algebras of compact operators}

In Part \ref{__2graded} (Sections \ref{_5tensor}-\ref{_8consequences}), we
consider general finitely subgraded Lie algebras of compact operators. Results
establishing the existence of invariant subspaces for these algebras are very
interesting, in the sense that hypotheses of a highly algebraic nature are
able to imply the existence of intrinsically geometric objects.

In Section \ref{_5tensor} a special tensor product construction is used to
reduce the consideration of $\Gamma$-subgraded Lie algebras of compact
operators to $\Gamma$-graded ones. More concretely, let $\pi$ be a faithful
representation of $\Gamma$ on a finite dimensional space $\mathcal{Y}$. To
every $\Gamma$-subgraded Lie algebra $L$ of operators on a Banach space
$\mathcal{X}$ one can associate a $\Gamma$-graded Lie algebra $L^{\pi}$ of
operators on $\mathcal{X}\otimes\mathcal{Y}$, with components $L_{\gamma
}\otimes\pi(\gamma)$ for $\gamma\in\Gamma$. We say that $L^{\pi}$ is a
\textit{graded ampliation} of $L$.

If $L$ consists of compact operators, then it turns out that $L$ is Engel or
$E$-solvable whenever $L^{\pi}$ is. This allows us to establish the result
that $L$ is Engel if the closure of each component of $L$ has no nonzero
finite rank operators. For cyclic $\Gamma$, we prove that $L$ is $E$-solvable
if the closure of every component contains no nonzero nilpotent finite rank
operators. This is a direct extension of the above-mentioned result from
\cite{ShT2005}.

Thus, in considering a subgraded Lie algebra of compact operators, we may
restrict ourselves to the case when Lie algebra contains nonzero homogeneous
finite rank operators. Such operators generate an ideal which also a subgraded
Lie algebra. In Section \ref{_7triangularization} we reduce the study of
subgraded Lie algebras of finite rank operators to the case when the
underlying space is finite dimensional. In Section \ref{_6finite}, we consider
this case in great detail.

Results on graded Lie algebras of operators on a finite dimensional space can
be considered as an aspect of the purely algebraic theory of graded groups and
Lie rings; for a treatment of the theory in this context, see \cite{Hig, Tho,
KrK, Kre, STW, KhM, MKh}. It was shown in \cite{KrK, Kre} that every
$\mathbb{Z}_{n}$-graded Lie ring $L$ with $L_{0}=0$ is solvable, with derived
length less than $2^{n-1}-1$, and with $L$ being nilpotent whenever $n$ is
prime. These estimates were improved in \cite[Theorems A, B and C]{STW}.

Here we do not require estimates on the derived length; for the convenience of
readers we provide, in Lemma \ref{prime}, a self-contained and comparatively
simple proof of the solvability of $L$ when $L$ consists of operators on a
finite dimensional space with scalar $L_{0}$. Another result, Lemma
\ref{finsubgraded}, is that a $\Gamma$-subgraded Lie algebra $L$ of operators
on a finite dimensional space is solvable if all its homogeneous elements are
Engel elements of $L$; in the case when $\Gamma=\mathbb{Z}_{n}$, it suffices
to assume that only $L_{0}$ consists of Engel elements of $L$. For
applications in Part \ref{__3related}, we require Lemma \ref{L0-triang}, that
if $L_{0}$ is solvable and noncommutative, then $L$ has a non-scalar solvable
ideal, and as a consequence, is reducible.

The main result of Part \ref{__2graded} is Theorem \ref{main1}; it states that
if $L$ is a $\Gamma$-subgraded Lie algebra of compact operators, and if the
homogeneous elements of $L$ are Volterra operators, or more generally, are
Engel elements of $L$, then $L$ is triangularizable. For $\Gamma
=\mathbb{Z}_{n}$, it is sufficient to require only that the component $L_{0}$
consists of Engel elements.

In Section \ref{_8consequences} we discuss some of the consequences of this
result. It is well known to algebraists that for any $\mathbb{Z}_{n}$-grading
of a Lie algebra $L$, one can find an automorphism $\phi$ of $L$ satisfying
the condition $\phi^{n} = 1$ such that the components of $L$ are the spectral
subspaces of $\phi$. Conversely, each $\phi$ defines a $\mathbb{Z}_{n}%
$-grading of $L$. Thus when $L$ is $\mathbb{Z}_{n}$-graded, the statement of
Theorem \ref{main1} can be reformulated in the following way (see Corollary
\ref{aut}(i)): if $\phi$ is an automorphism of finite order of a closed Lie
algebra $L$ of compact operators, and if every element invariant under $\phi$
is a Volterra operator, then $L$ is triangularizable. In Corollary
\ref{aut}(ii) we extend this result to automorphisms of possibly infinite
order, when the spectrum of $\phi$ be contained in a finite subgroup of the
unit circle. Finally, Corollary \ref{subgroup} combines both statements of the
main result in the following way: if $\Gamma$ has a subgroup $H$ such that
$\Gamma/H$ is cyclic, and such that every elements in $L_{\gamma}$ is Engel,
for $\gamma\in H$, then $L$ is triangularizable.

It should be noted that one can obtain further results by imposing
quasinilpotence conditions on certain other ``valuable'' Lie multiplicative
sets of homogeneous elements. The first result of this kind is Theorem
\ref{multiset}, which states that for the triangularizability of a $\Gamma
$-subgraded Lie algebra $L$ of compact operators, it suffices to require that
all the homogeneous commutators are Engel elements of $L$. We next prove a
stronger assertion (see Corollary \ref{main2}) in which the quasinilpotence
condition is imposed only for homogeneous commutators in certain components
(depending on $\Gamma$). For instance, if $\Gamma=\mathbb{Z}_{n}$, it suffices
to require that every commutator $[a,b]$ is an Engel element of $L$, for $a\in
L_{k}$, $b\in L_{-k}$, and $k\in\mathbb{Z}_{n}$. Of course, these conditions
are also necessary for the triangularizability of $L$.

An easy consequence of these results on homogeneous commutators is that for a
$\Gamma$-subgraded Lie algebra $L$ of compact operators, homogeneous local
triangularizability implies global triangularizability. The simplest form of
this statement, when every pair of homogeneous elements is triangularizable,
is shown in Corollary \ref{loctriang}. A stronger result, requiring only the
homogeneous local triangularizability of certain components (again, depending
on $\Gamma$), is given in Corollary \ref{tria}.

At the end of Part \ref{__2graded} we discuss an extension of our results to
the case when $\Gamma$ is an infinite commutative group (see Theorem
\ref{locfin}).

\subsubsection{$\mathbb{Z}_{2}$-subgraded Lie algebras, Lie triple systems,
and Jordan algebras}

In Part \ref{__3related}, $\mathbb{Z}_{2}$-subgraded Lie algebras, Lie triple
systems and Jordan algebras are treated. As was mentioned above, in the case
of $\mathbb{Z}_{2}$-subgraded Lie algebra $L=L_{0}+L_{1}$, we consider the
case when the condition of quasinilpotence is imposed on the component $L_{1}%
$, and not on $L_{0}$. (The case when the condition of quasinilpotence is
imposed on $L_{0}$ was treated in Part \ref{__2graded}.) Note that
$[L_{1},[L_{1},L_{1}]]\subset L_{1}$, which means that $L_{1}$ is a Lie triple system.

The main result of Part \ref{__3related}, Theorem \ref{two}, states that if
$L$ is a $\mathbb{Z}_{2}$-subgraded Lie algebra of compact operators, and if
$L_{1}$ is non-scalar and consists of Engel elements of $L$ (in particular, if
$L_{1}$ consists of Volterra operators), then $L$ is reducible. Of course, if
in addition $L_{1}$ generates $L$ as a Lie algebra then $L$ is
triangularizable. Example \ref{e2} shows that this is not valid for $n=4$,
i.e. there exists an irreducible $\mathbb{Z}_{4}$-subgraded Lie algebra of
finite rank operators such that every element in $L_{1}$ is nilpotent, and
$L_{1}$ generates $L$ as a Lie algebra. It seems that the corresponding
problem for $n=3$ is still open.

Applying the results of Part \ref{__2graded} to the case when $L$ is a
$\mathbb{Z}_{2}$-subgraded Lie algebra of compact operators, we get in Theorem
\ref{2theor}(ii) that if $[L_{0},L_{0}]$ and $[L_{1},L_{1}]$ consist of
Volterra operators, then $L$ is triangularizable.

Every Lie algebra or Jordan algebra of operators forms a Lie triple system. We
investigate operator Lie triple systems, and prove the result that a Lie
triple system of Volterra operators is triangularizable (see Theorem
\ref{tripvolt}). This is an infinite dimensional extension of a result in
\cite{Ja8}.

The proof of this result for operator Lie triple systems illustrates the
benefits of working with subgraded, and hence not necessarily graded, Lie
algebras. In fact, if $L_{1}$ is an operator Lie triple system then
$L_{1}+\left[  L_{1},L_{1}\right]  $ is a $\mathbb{Z}_{2}$-subgraded Lie
algebra when $L_{0}=\left[  L_{1},L_{1}\right]  $ (see Lemma \ref{Jor-grad}),
but it is possible that $L_{0}\cap L_{1}\neq0$.

As an application, we establish that a Jordan algebra of Volterra operators is
triangularizable (see Theorem \ref{jordvolt}), and moreover, that a Jordan
algebra of compact operators which has a nonzero Jordan ideal of Volterra
operators is reducible (see Theorem \ref{jordan}). These results are the
Jordan analogues of results for Lie algebras which were established in
\cite{ShT2000} and \cite{ShT2005} respectively.

The authors would like to express their gratitude to Matej Bre\v{s}ar and
Heydar Radjavi, who provided them with helpful discussions and moral support.

\part{\label{__1criteria}Triangularization Criteria for Lie Algebras of
Compact Operators}

In Theorem \ref{redcri} we collect the known necessary and sufficient
conditions for the reducibility and triangularizability of Lie algebras of
compact operators. Our goal here is to add to these conditions in order to
prove our main results. This is done in Theorem \ref{crit} and Remark
\ref{addcrit}.

Many other results of this part can be treated as preliminary ones for
consideration of subgraded Lie algebras of compact operators.

\section{\label{_2preliminaries} Preliminaries}

For a complex Banach space $\mathcal{X}$, we denote by $\mathcal{B}%
(\mathcal{X})$ the algebra of all bounded linear operators on $\mathcal{X}$,
and by $\mathcal{K}(\mathcal{X})$ (respectively, $\mathcal{F}(\mathcal{X})$)
the ideal of all compact (respectively, finite rank) operators on
$\mathcal{X}$. Recall that an operator is called \textit{Volterra} if it is
compact and quasinilpotent. A set of operators is called \textit{Volterra} if
all its elements are Volterra.

We will implicitly use the well known fact (see, for instance, \cite{Mu1})
that the spectrum is continuous on $\mathcal{K}(\mathcal{X})$, and more
generally, on the set of operators with countable spectra. In particular, the
limit of a convergent sequence of Volterra operators is Volterra.

A\ set $M$ of bounded linear operators on a Banach space $\mathcal{X}$ is said
to be \textit{triangularizable} if there exists a maximal subspace chain
consisting of closed subspaces that are invariant under all operators in $M$.
The usage of the word triangularizable stems from the fact that if
$\mathcal{X}$ is finite dimensional, the existence of such a subspace chain is
equivalent to the existence of basis of $\mathcal{X}$ under which all
operators in $M$ can be represented by upper-triangular matrices. As with
matrices, it is often much easier to work with triangularizable families of
operators, and in fact, it could be said that one of the aims of Invariant
Subspace Theory is to establish triangularizability results.

It is easy to see that if a property of a set $M$ of operators implies
reducibility of $M$ and is preserved for the related set $M|\mathcal{V}$ of
induced operators on the quotient $\mathcal{V}=\mathcal{Y}/\mathcal{Z}$ of
arbitrary invariant closed subspaces $\mathcal{Y}$ and $\mathcal{Z}$ for $M$,
with $\mathcal{Z}\subset\mathcal{Y}$, then this property implies
triangularizability of $M$ (this is the Triangularization Lemma in
\cite{RR0}). That is why to establish the triangularizability of a family of
operators in many cases it actually suffices to establish its reducibility.

For a set $M$ of compact operators, \textit{the triangularizability of} $M$
\textit{is equivalent to the commutativity of} $\mathcal{A}(M)$ \textit{modulo
the Jacobson radical} $\operatorname*{rad}\mathcal{A}(M)$ (see \cite{Mu2}),
where $\mathcal{A}(M)$ is the closed algebra generated by $M$. This implies
important spectral properties of $M$, for instance that the spectral radius
(and also the spectrum itself) is subadditive and submultiplicative on
$\mathcal{A}(M)$, that commutators of elements of $\mathcal{A}(M)$ are
Volterra, etc. In what follows, we will frequently make use of these properties.

A linear subspace $L$ of $\mathcal{B}(\mathcal{X})$ is called an
\textit{operator Lie algebra} if it is closed under the Lie product
\[
\lbrack a,b]=ab-ba.
\]
For each $a\in L$, we denote by $\operatorname*{ad}a$, or $\operatorname*{ad}%
_{L}a$, the operator on $L$ which maps $x$ to $[a,x]$, for every $x\in L$.

For brevity, we will occasionally make use of shorthand notation. For example,
we will write $UV$ to denote the (algebraic) associative product of linear
subspaces $U$ and $V$ of $\mathcal{B}(\mathcal{X})$, i.e.%
\[
UV=\operatorname*{span}\{ab:a\in U,b\in V\},
\]
and similarly, we will write $[U,V]$ to denote the (algebraic) Lie product,
i.e.
\[
\left[  U,V\right]  =\operatorname*{span}\left\{  \left[  a,b\right]  :a\in
U,b\in V\right\}  .
\]

Recall that a finite dimensional Lie algebra $L$ is called \textit{semisimple}
if it does not have non-zero abelian ideals. Cartan's well-known
Semisimplicity Criterion (see for example \cite[Section 3.4]{Jac}) states that
$L$ \textit{is semisimple if and only if the Killing form }%
\[
\left\langle x,y\right\rangle :=\operatorname*{tr}((\operatorname*{ad}%
x)(\operatorname*{ad}y))
\]
\textit{ of }$L$ \textit{is non-degenerate on }$L$; here $\operatorname*{tr}%
(a)$ means the trace of operator $a$.

For a Lie algebra $L$, the \textit{lower central} (respectively,
\textit{derived}) \textit{series} of ideals $L^{\left[  k\right]  }$
(respectively, $L^{(k)}$) is defined \cite[Section 1.3]{Hum} by the rules
\begin{align*}
L^{\left[  0\right]  }  &  =L\quad\text{and\quad}L^{\left[  k+1\right]
}=[L,L^{\left[  k\right]  }]\text{ }\\
\text{(respectively, }L^{\left(  0\right)  }  &  =L\quad\text{and}\quad
L^{(k+1)}=[L^{(k)},L^{(k)}]\text{)}%
\end{align*}
for $k=0,1,\ldots$. Note that we use the notation $L^{\left[  k\right]  }$
instead of the standard $L^{k}$ to differentiate ideals of the lower central
series from corresponding associative structures.

If $L^{\left[  k\right]  }=0$ (respectively, $L^{(k)}=0$) for some $k$, then
$L$ is called \textit{nilpotent} (respectively, \textit{solvable}). One of the
forms of Lie's classical Theorem is that \textit{a Lie algebra of operators on
a finite-dimensional space is solvable if and only if it is triangularizable
}(or more precisely, is representable by upper-triangular matrices in some
basis of the space).

The following lemma of Cartan, which is also well-known, provides a criterion
for the solvability of $L$ (see \cite[Theorem 1.5.2]{Bou} and \cite[Section
4.3]{Hum}).

\begin{lemma}
[Cartan's Solvability Criterion]\label{cartan}Let $L$ be a Lie algebra of
operators on a finite dimensional vector space. Then the following conditions
are equivalent.

\begin{itemize}
\item $L$ is solvable.

\item $\left\langle a,b\right\rangle =0$ for all $a\in\lbrack L,L]$ and $b\in
L$.

\item $\operatorname*{tr}(ab)=0$ for all $a\in\lbrack L,L]$ and $b\in L$.
\end{itemize}
\end{lemma}

In particular, it follows from Cartan's Solvability Criterion that \textit{if
}$\operatorname*{tr}(ab)=0$\textit{ for all }$a$\textit{,}$b\in L$\textit{
then }$[L,L]$\textit{ consists of nilpotent operators.} This was extended in
\cite{Ken} to infinite dimensions as follows.

\begin{theorem}
\cite[Theorem 4.6]{Ken}\label{carMatt} Let $L$ be a Lie algebra of trace class
operators (on a Hilbert space) satisfying $\operatorname*{tr}(ab)=0$ for all
$a,b\in L$. Then every finite rank operator in $[L,L]$ is nilpotent.
\end{theorem}

For subspaces $N$ and $M$ of $\mathcal{B}(\mathcal{X})$, we will write
\begin{equation}
\operatorname*{tr}(NM)=0 \label{1-trace}%
\end{equation}
to indicate that $\operatorname*{tr}(ab)=0$ holds for every $a\in N$ and $b\in
M$, requiring of course that each $ab$ is an operator with trace. Note that in
virtue of linearity of trace, condition (\ref{1-trace}) is equivalent to the
condition $\operatorname*{tr}(c)=0$ for all $c\in NM$.

Similarly, we will write
\[
\left\langle N,M\right\rangle =0
\]
to indicate that $\left\langle a,b\right\rangle =0$ for every $a\in N$ and
$b\in M$.

Let us give a somewhat more general formulation of Theorem \ref{carMatt}. It
should be noted that its proof doesn't need any modification for operators in
a normed operator ideal (see \cite{Pie}) with \textit{spectral trace }(recall
that the trace is \textit{spectral} if it coincides with the sum of
eigenvalues, taken with multiplicity).

\begin{lemma}
\label{carfin} Let $L$ be a Lie algebra of compact operators. If there is a
non-zero finite rank operator $a$ in $L$ such that $\operatorname*{tr}(aL)=0$,
then $L$ contains either a non-scalar commutative ideal, or a non-zero ideal
consisting of nilpotent operators. In both cases, $L$ is reducible.
\end{lemma}

\begin{proof}
Let
\[
I=\{x\in L\cap\mathcal{F}(\mathcal{X}):\operatorname*{tr}(xL)=0\}.
\]
Then $I$ is an ideal of $L$. Indeed, for arbitrary $x\in I$ and $y\in L$, we
have that%
\[
\operatorname*{tr}([x,y]L)=\operatorname*{tr}(x[y,L])=0\text{.}%
\]
Therefore, we obtain that
\[
\operatorname*{tr}(II)=0.
\]
It follows from Theorem \ref{carMatt} that $[I,I]$ consists of nilpotent
operators. If $[I,I]\not =0$, then it is a non-zero Volterra ideal of $L$.
Otherwise, $[I,I]=0$, and we have that $I$ is a commutative ideal of $L$
containing a non-scalar operator $a$ (since $\operatorname*{tr}(a^{2})=0$). In
both cases, $L$ is reducible by \cite[Theorem 4.26]{ShT2005} (see Theorem
\ref{solid} below).
\end{proof}

For an operator $a\in\mathcal{B}(\mathcal{X})$ and $\lambda\in\mathbb{C}$,
let
\[
\mathcal{E}_{\lambda}(a)=\{x\in\mathcal{X}:\lim\left\Vert \left(
a-\lambda\right)  ^{n}x\right\Vert ^{1/n}=0\}.
\]
This linear subspace of $\mathcal{X}$ is an\textbf{ }\textit{elementary
spectral manifold }of $a$ \cite[Section 3]{ShT2005}. It is a useful fact (see
for instance \cite[Proposition 3.3]{ShT2005}) that $\mathcal{E}_{\lambda}(a)$
is always closed if $a$ has countable spectrum (or, more generally, if $a$ is
decomposable). Thus, if a decomposable operator $a$ is \textit{locally
quasinilpotent }on a subset $G$ of $\mathcal{X}$ (this means that
\[
\lim\left\Vert a^{n}x\right\Vert ^{1/n}=0
\]
for every $x\in G$), $a$ is locally quasinilpotent, and hence is
quasinilpotent, on the closed linear span of $G$ in virtue of well-known fact
that \textit{if a bounded operator is locally quasinilpotent on a Banach space
then it is quasinilpotent. }We use in the paper a stronger version of this
result (see for instance \cite[Corollary 14.19]{M}), namely the spectral
radius $\rho(a)$ of an operator $a\in\mathcal{B}(\mathcal{X})$ satisfies to
the equality
\begin{equation}
\rho(a)=\underset{x\in\mathcal{X}}{\sup}\lim\sup\left\Vert a^{n}x\right\Vert
^{1/n}.\label{f_muller}%
\end{equation}

Recall that an element $a$ of a normed Lie algebra $L$ is called an
\textit{Engel element} of $L$ if the operator $\operatorname*{ad}a$ is
quasinilpotent on $L$, i.e. if
\[
\lim\left\Vert (\operatorname*{ad}a)^{n}\right\Vert ^{1/n}=0.
\]
It is clear that if $L$ is an operator Lie algebra, and if $a$ is
quasinilpotent, then $a$ is an Engel element of $L$, but the converse is not
true in general. We say that $L$ is \textit{Engel} if all its elements are Engel.

If $a$ is an element of a normed Lie algebra $L$, the following useful
property
\[
\lbrack\mathcal{E}_{\lambda}(\operatorname*{ad}a),\mathcal{E}_{\mu
}(\operatorname*{ad}a)]\subset\mathcal{E}_{\lambda+\mu}(\operatorname*{ad}a)
\]
holds for every $\lambda,\mu\in\mathbb{C}$; in the case when  $L$ is a normed
algebra we have also that
\begin{equation}
\mathcal{E}_{\lambda}(\operatorname*{ad}a)\mathcal{E}_{\mu}(\operatorname*{ad}%
a)\subset\mathcal{E}_{\lambda+\mu}(\operatorname*{ad}a)\label{f_addition}%
\end{equation}
(see a more general fact in \cite[Lemma 3.5]{ShT2005}). In what follows we
will also need the fact that\textit{ if} $a$ \textit{is a compact operator,
and} $\lambda\neq0$, \textit{then }$\mathcal{E}_{\lambda}(\operatorname*{ad}%
a)$ \textit{consists of nilpotent finite rank operators; } the inclusion
\[
\mathcal{E}_{\lambda}(\operatorname*{ad}a)\subset\mathcal{F}(\mathcal{X})
\]
for $\lambda\neq0$ was proved in \cite[Theorem 3]{Woj}, and we refer the
nilpotency of elements of $\mathcal{E}_{\lambda}(\operatorname*{ad}a)$ for
$\lambda\neq0$ to mathematical folklore. For instance, it is an easy
consequence of (\ref{f_addition}) and boundedness of $\sigma
(\operatorname*{ad}a)$. It is convenient to reformulate this fact as follows
(see, for example, \cite[Corollary 2]{ShT1999}).

\begin{theorem}
\label{engel} Let $L$ be a closed Lie subalgebra of compact operators on an
infinite dimensional space $\mathcal{X}$. If $L$ doesn't contain non-zero
nilpotent finite rank operators, then $L$ is Engel.
\end{theorem}

Recall that a normed Lie algebra is called \textit{Engel-solvable,} or
$E$\textit{-solvable,} if, for each proper closed ideal $I$ of $L$, the
quotient Lie algebra $L/I$ contains a non-zero closed Engel ideal; an ideal of
a normed Lie algebra is called \textit{Engel }(respectively, $E$%
\textit{-solvable}) if it is an Engel (respectively, $E$-solvable) Lie algebra.

It was proved in \cite[Corollary 4.25]{ShT2005} that, for a Lie algebra $L$ of
compact operators, the condition of $E$-solvability is equivalent to the
geometric condition of triangularizability. As we will use this result
frequently in what follows, it is convenient to collect in one statement
several known conditions in \cite{ShT2005} that are equivalent to $L$ being
$E$-solvable; note that these conditions and proofs of their implications are
distributed in \cite[Corollaries 4.20, 4.25, 5.13 and 5.16, Theorem 5.19,
etc]{ShT2005}. To aid the reader, we spell out exactly the nontrivial
implications contained in the following theorem: $(1)\iff(2)$ is a part of
\cite[Corollary 4.25]{ShT2005}, $(3)\Longrightarrow(2)$ follows immediately
from \cite[Corollary 5.16 and Theorem 5.19]{ShT2005}, $(5)\iff(2)\iff(6)$ is
an easy consequence of \cite[Theorems 5.2 and 5.9, Corollaries 5.13 and
5.16]{ShT2005}. All other implications required, namely $(2)\Longrightarrow
(4)\Longrightarrow(3)$, are trivial.

\begin{theorem}
\label{redcri} For a Lie subalgebra $L$ of $\mathcal{K}(\mathcal{X})$ the
following conditions are equivalent.

\begin{enumerate}
\item $L$ is $E$-solvable.

\item \label{redcri1} $L$ is triangularizable.

\item $[L,L]$ is Engel.

\item $[L,L]$ is Volterra.

\item \label{redcri5}The spectral radius is subadditive on $\operatorname*{ad}%
L$.

\item \label{redcri6}The spectral radius is subadditive on $L$ (that is,
$\rho\left(  a+b\right)  \leq\rho\left(  a\right)  +\rho\left(  b\right)  $
for every $a,b\in L$).
\end{enumerate}
\end{theorem}

An immediate consequence of the above results is that an Engel Lie algebra of
compact operators is triangularizable.

The following result will be our main tool for establishing reducibility results.

\begin{theorem}
\label{solid} \cite[Theorem 4.26]{ShT2005} Let $L$ be a Lie subalgebra of
$\mathcal{K}(\mathcal{X})$. If there is a non-scalar $E$-solvable ideal of $L$
then $L$ is reducible.
\end{theorem}

Here, a family of operators is \textit{non-scalar }if it contains a non-scalar
operator (recall that operator is \textit{scalar }if it is a scalar multiple
of the identity operator). Trivially, \textit{a non-zero compact operator on
}$\mathcal{X}$\textit{ is non-scalar if }$\mathcal{X}$\textit{ is infinite
dimensional.}

We will also need the following structural result.

\begin{theorem}
\label{twoideals}\cite[in Theorem 4.32]{ShT2005} Let $L$ be a Lie subalgebra
of $\mathcal{K}(\mathcal{X})$. If $L$ has two non-scalar ideals with scalar
intersection then $L$ is reducible.
\end{theorem}

\textit{For finite-dimensional Lie algebras}, one can easily check that $L$
\textit{is Engel} (respectively, $E$\textit{-solvable}) \textit{if and only if
it is nilpotent} (respectively,\textit{ solvable}).

It is well known (see for example \cite{Jac}) that any finite-dimensional Lie
algebra $L$ has a largest solvable ideal $R$ (the \textit{radical} of $L$),
and that the quotient Lie algebra $L/R$ is semisimple. A similar result for
existence of the largest $E$-solvable ideal was proved in \cite[Corollary
5.15]{ShT2005} for a Lie algebra of compact operators. We note the following
important property of this ideal.

\begin{theorem}
\label{inner}\cite[Corollary 5.15]{ShT2005} Let $L$ be a Lie algebra of
compact operators, and let $R$ be the largest $E$-solvable ideal of $L$. Then
$R$ is inner-characteristic, i.e. $[a,R]\subset R$ for every $a\in
\mathcal{B}(\mathcal{X})$ such that $[a,L]\subset L$.
\end{theorem}

\section{\label{_3subalgebra}Engel Elements of a Lie Algebra in its
Subalgebra}

In this section, we will consider a Lie algebra $L$ of compact operators
having a subalgebra $L_{0}$. We will be particularly interested in Engel
elements of $L$ that also belong to $L_{0}$.

In applications $L$ will be subgraded and $L_{0}$ will be its $0$-component.

\begin{lemma}
\label{space}Let $L$ be a closed Lie algebra of compact operators, $L_{0}$ a
closed subalgebra of $L$, $N$ the set of all elements of $L_{0}$ that are
Engel elements of $L$ (or the set of all Volterra operators in $L_{0}$),
$N_{f}$ the set of all finite-rank operators in $N$. If $N$ (respectively
$N_{f}$) is a subspace of $L_{0}$ then it is an ideal of $L_{0}$.
\end{lemma}

\begin{proof}
Let $b\in L_{0}$. Then $\exp\left(  \lambda\operatorname*{ad}b\right)  $ is a
bounded automorphism of $L$ for every $\lambda\in\mathbb{C}$. If $a$ is in $N$
then
\[
\left(  \exp\left(  \lambda\operatorname*{ad}b\right)  \right)  (a)=\exp
(\lambda b)a\exp(-\lambda b)\in N
\]
for every $\lambda\in\mathbb{C}$, whence we obtain by hypothesis that
\[
\upsilon(\lambda):=\lambda^{-1}\left[  \left(  \exp\left(  \lambda
\operatorname*{ad}b\right)  \right)  (a)-a\right]
\]
is in $N$ for every $\lambda\neq0$, and that
\[
\left[  b,a\right]  =\lim_{\lambda\rightarrow0}\upsilon(\lambda)
\]
is in $N$. Note that if $a\in N_{f}$ then so is $\left[  b,a\right]  $.
\end{proof}

We will need the following three lemmas about Banach algebras.

\begin{lemma}
\label{sum}Let $A,B$ be Banach algebras, $\phi:A\longrightarrow B$ a bounded
homomorphism. Suppose that $A$ consists of elements with countable spectra,
and that $A/\operatorname*{rad}A$ is commutative. If $\phi(x)$ and $\phi(y)$
are quasinilpotent elements of $B$ for $x,y\in A$ then $\phi(x+y)$ is
quasinilpotent, too.
\end{lemma}

\begin{proof}
One may assume that $A$ and $B$ are unital and $\phi(1)=1$. It is clear that
$A/\ker\phi$ is commutative modulo the Jacobson radical and consists of
elements with countable spectra. So we can assume that $\phi$ is injective.

We claim that $x,y$ \textit{are quasinilpotent elements of} $A$.

Indeed, if $x$ is not quasinilpotent then there is a non-zero isolated point
$\lambda\in\sigma(x)$. Let $p$ be the Riesz idempotent of $x$ corresponding to
$\lambda$. Then
\[
p\neq0
\]
and $(x-\lambda)p$ is quasinilpotent, whence so is $\left(  \phi
(x)-\lambda\right)  \phi(p)$. If $\phi(p)\neq0$ then $\lambda\in\sigma
(\phi(x))$, a contradiction. Therefore
\[
\phi(p)=0.
\]
But, since $\phi$ is injective, $p=0$, a contradiction. This proves our claim.

Since $A$ is commutative modulo the Jacobson radical, $x+y$ is a
quasinilpotent element of $A$. Then $\phi(x+y)$ is quasinilpotent, too.
\end{proof}

Recall \cite{Ale} that a Banach algebra $A$ is called \textit{compact} if the
map \textrm{L}$_{x}\mathrm{R}_{x}:y\longmapsto xyx$ is compact on $A$, for
every $x\in A$, where multiplication operators $\mathrm{L}_{x}$ and
\textrm{R}$_{x}$ in $\mathcal{B}(A)$ are defined by
\[
\mathrm{L}_{x}y=xy\quad\text{and\quad}\mathrm{R}_{x}y=yx
\]
for every $y\in A$. It should be noted \cite[Theorem 4.4]{Ale} that each
element of a compact Banach algebra has countable spectrum, and every closed
subalgebra of $\mathcal{K}(\mathcal{X})$ is a compact Banach algebra.

The following result is contained in \cite[Theorems 2 and 7]{ShT2001}. Since
the proofs of these statements are still not published, we include the proof
adapted for this special case.

\begin{lemma}
\label{tens} Let $A$ and $B$ be compact Banach algebras, and let
$A^{1}=A\oplus\mathbb{C}$ and $B^{1}=B\oplus\mathbb{C}$ be the algebras
obtained by adjoining the identity element. If $A$ and $B$ are commutative
modulo the Jacobson radical, then so is the projective tensor product
$W=A^{1}\widehat{\otimes}B^{1}$, and all elements of $W$ have countable spectra.
\end{lemma}

\begin{proof}
Let $I=A\widehat{\otimes}B$. As is known, under the natural embedding of $I$
into $W$, the algebra $I$ is topologically isomorphic to a closed subspace of
$W$ which is an ideal of $W$. So we will identify $I$ with a closed ideal of
$W$.

Let $\tau$ be a strictly irreducible representation of $W$. If $\tau(I)=0$
then $\tau$ is also a strictly irreducible representation of the algebra
$\left(  A\oplus B\right)  ^{1}$ (with coordinate-wise multiplication). This
algebra is clearly commutative modulo the Jacobson radical, whence $\tau$ is one-dimensional.

Thus, one can assume that $\tau(I)\neq0$. Then $\tau$ is also a strictly
irreducible representation of $I$. Using the fact that the tensor product of
compact operators is a compact operator (on the projective tensor product of
corresponding spaces), we obtain that $I$ is a compact Banach algebra. By
\cite[Theorem 5.1]{Ale}, $\tau(I)$ consists of compact operators. Clearly, the
set
\[
G=\left\{  x\otimes y:x\in\operatorname*{rad}A,y\in B\right\}  \cup\left\{
z\otimes w:z\in A,w\in\operatorname*{rad}B\right\}
\]
is a multiplicative semigroup of quasinilpotent elements of $I$. Therefore
$\tau(G)$ is a Volterra semigroup. Let $N$ be the closure of
$\operatorname*{span}\tau(G)$ in $\tau(I)$. By \cite[Theorem 4]{Tur1998}, $N$
is Volterra. But $N$ is an ideal of $\tau(I)$, whence
\[
N=0.
\]
Since it is easy to check that, for every $x,z\in A$ and $y,w\in B$,%
\[
\left[  x\otimes y,z\otimes w\right]  =xz\otimes yw-zx\otimes wy=\left[
x,z\right]  \otimes yw+zx\otimes\left[  y,w\right]
\]
with $[x,z]\in\operatorname*{rad}A$ and $\left[  y,w\right]  \in
\operatorname*{rad}B$, we obtain that
\[
\tau(\left[  I,I\right]  )\subset N.
\]
Then $\tau(I)$ is commutative and $\tau$ is one-dimensional.

In any case, we obtain that $W$ is commutative modulo the Jacobson radical.

Now we have to prove that $W$ consists of elements with countable spectra. For
every $g\in W$ there are $x\in A$, $y\in B$, $\lambda\in\mathbb{C}$, and $f\in
I$ such that%
\[
g=\lambda(1\otimes1)+x\otimes1+1\otimes y+f.
\]
Then, by commutativity of $W$ modulo the Jacobson radical, we obtain that
\[
\sigma(g)\subset\lambda+\sigma(x\otimes1)+\sigma(1\otimes y)+\sigma
(f)\subset\lambda+\sigma(x)+\sigma(y)+\sigma(f)
\]
Since all last spectra are countable (in particular, $f$ has countable
spectrum by \cite[Theorem 4.4]{Ale}), $g$ has countable spectrum.
\end{proof}

Recall that an element $b$ of a Banach algebra $B$ is\textit{ of finite rank
}if the map \textrm{L}$_{b}\mathrm{R}_{b}:x\mapsto bxb$ is of finite rank on
$B$. By \cite[Theorem 7.2]{Ale}, for a semisimple compact Banach algebra $B$,
an element $b$ of $B$ is of finite rank if and only if $b$ is in the socle
$\operatorname*{soc}B$ of $B$. Of course, every quasinilpotent finite rank
element of $B$ is in fact nilpotent.

The following lemma strongly generalizes \cite[Theorem 3]{Woj}.

\begin{lemma}
\label{socle} Let $B$ be a semisimple compact Banach algebra $B$, and let
$a\in B$. Then $\mathcal{E}_{1}\left(  \operatorname*{ad}_{B}a\right)  $
consists of nilpotent elements of $\operatorname*{soc}B$.
\end{lemma}

\begin{proof}
Recall that the spectrum of $a$ is considered with respect to $B^{1}$. By
\cite[Theorem 4.4]{Ale}, there is a finite number of points in $\sigma
(a)\cap\mathbb{C}_{\varepsilon}$ for every $\varepsilon>0$, where
$\mathbb{C}_{\varepsilon}=\left\{  \mu\in\mathbb{C}:\left\vert \mu\right\vert
>\varepsilon\right\}  $, and by \cite[Lemma 3]{ShT2002}, the Riesz idempotent
of $a$ corresponding to a non-zero number $\lambda$ in $\sigma(a)$ is of
finite rank in $B$. Thus, if $p_{\varepsilon}$ is the Riesz idempotent of $a$
corresponding to $\sigma(a)\cap\mathbb{C}_{\varepsilon}$, then $p_{\varepsilon
}$ is of finite rank in $B$. Let $q_{\varepsilon}=1-p_{\varepsilon}$. Since
$a$ commutes with $q_{\varepsilon}$ and $\mathrm{L}_{q_{\varepsilon}%
}\mathrm{R}_{q_{\varepsilon}}$ is a projection, it is easy to verify that
\begin{equation}
\mathrm{L}_{q_{\varepsilon}}\mathrm{R}_{q_{\varepsilon}}\mathcal{E}_{1}\left(
\operatorname*{ad}a\right)  \subset\mathcal{E}_{1}\left(  \operatorname*{ad}%
a\right)  \cap\mathcal{E}_{1}\left(  \mathrm{L}_{q_{\varepsilon}}%
\mathrm{R}_{q_{\varepsilon}}\operatorname*{ad}a\right)  , \label{proj}%
\end{equation}
and since $\rho(aq_{\varepsilon})\leq\varepsilon$, that
\begin{align*}
\rho\left(  \mathrm{L}_{q_{\varepsilon}}\mathrm{R}_{q_{\varepsilon}%
}\operatorname*{ad}a\right)   &  =\rho\left(  \mathrm{L}_{aq_{\varepsilon}%
}\mathrm{R}_{q_{\varepsilon}}-\mathrm{L}_{q_{\varepsilon}}\mathrm{R}%
_{aq_{\varepsilon}}\right) \\
&  \leq\rho\left(  aq_{\varepsilon}\right)  \rho\left(  q_{\varepsilon
}\right)  +\rho\left(  q_{\varepsilon}\right)  \rho\left(  aq_{\varepsilon
}\right)  \leq2\rho\left(  aq_{\varepsilon}\right)  \leq2\varepsilon.
\end{align*}
If $\varepsilon<1/2$, then $\rho\left(  \mathrm{L}_{q_{\varepsilon}}%
\mathrm{R}_{q_{\varepsilon}}\operatorname*{ad}a\right)  <1$, and then
\[
\mathcal{E}_{1}\left(  \mathrm{L}_{q_{\varepsilon}}\mathrm{R}_{q_{\varepsilon
}}\operatorname*{ad}a\right)  =0,
\]
whence, by (\ref{proj}),
\[
\left(  1-p_{\varepsilon}\right)  b\left(  1-p_{\varepsilon}\right)  =0
\]
for every $b\in\mathcal{E}_{1}\left(  \operatorname*{ad}a\right)  $, that
completes the proof that $\mathcal{E}_{1}\left(  \operatorname*{ad}a\right)  $
is in $\operatorname*{soc}B$.

Since
\[
\left(  \mathcal{E}_{1}\left(  \operatorname*{ad}a\right)  \right)
^{n}\subset\mathcal{E}_{n}\left(  \operatorname*{ad}a\right)
\]
for every integer $n>0$ by \cite[Lemma 3.5]{ShT2005}, and since $\sigma\left(
\operatorname*{ad}a\right)  $ is bounded, we may take an integer $n$ outside
of $\sigma\left(  \operatorname*{ad}a\right)  $ and obtain that%
\[
\mathcal{E}_{n}\left(  \operatorname*{ad}a\right)  =0,
\]
whence every element of $\mathcal{E}_{1}\left(  \operatorname*{ad}a\right)  $
is nilpotent.
\end{proof}

For an algebra $B$, we denote by $B^{\operatorname*{op}}$ the algebra that is
opposite to $B$; this algebra coincides with $B$ as a linear space, but has
opposite multiplication, namely
\[
a\cdot b=ba
\]
for every $a,b\in B^{\operatorname*{op}}$. It is clear from this that if $B$
is a compact Banach algebra then so is $B^{\operatorname*{op}}$.

\begin{theorem}
\label{tensor}Let $L$ be a Lie algebra of compact operators, $L_{0}$ an
$E$-solvable subalgebra of $L$. If $a,b\in L_{0}$ are Engel elements of $L$
then so is $a+b$.
\end{theorem}

\begin{proof}
Let $D=\mathcal{A}(L_{0}\cup\{1\})$ and $A=D\widehat{\otimes}%
D^{\operatorname*{op}}$. Since $L_{0}$ is triangularizable by Theorem
\ref{redcri}, then the algebras $D$ and $D^{\operatorname*{op}}$ are
commutative modulo the Jacobson radical by \cite{Mu2}. Since $\mathcal{A}%
(L_{0})$ is a compact Banach algebra by \cite{Val} and \cite[Lemma 3.4]{Ale},
then $A$ is commutative modulo the Jacobson radical, and elements of $A$ have
countable spectra, by Lemma \ref{tens}.

Let $C=\mathcal{A}(L)$ and $B=\mathcal{B}(C)$. For every $x\in C$, it is
convenient to work with left and right multiplication operators $\mathrm{L}%
_{x}$ and \textrm{R}$_{x}$ on $C$. Define a map $\phi:A\longrightarrow B$ by
setting%
\[
\phi(x\otimes y)=\mathrm{L}_{x}\mathrm{R}_{y}%
\]
for every $x\otimes y\in A$, and by extending the map to the whole $A$ by
linearity and continuity. Then $\phi$ is a bounded homomorphism of Banach
algebras and $\phi(1\otimes1)$ is the identity operator on $C$.

We claim that \textit{if} $a$ \textit{is an Engel element of} $L$
\textit{then} $a$ \textit{is an Engel element of} $\mathcal{A}(L)$. Indeed,
since $\mathcal{E}_{0}\left(  \operatorname*{ad}_{C}a\right)  $ is a closed
algebra (see, for instance, \cite[Proposition 3.3 and Lemma 3.5]{ShT2005})
and
\[
L\subset\mathcal{E}_{0}\left(  \operatorname*{ad}\nolimits_{L}a\right)
\subset\mathcal{E}_{0}\left(  \operatorname*{ad}\nolimits_{C}a\right)  ,
\]
$\mathcal{E}_{0}\left(  \operatorname*{ad}\nolimits_{C}a\right)  $ contains
the closed algebra generated by $L$. This means that $a$ is an Engel element
of $\mathcal{A}(L)$.

Now let $a,b\in L_{0}$ be Engel elements of $L$. Since $C=\mathcal{A}(L)$,
they are Engel elements of $C$ by above. Take elements
\[
x=a\otimes1-1\otimes a\text{\quad and\quad}y=b\otimes1-1\otimes b
\]
in $A$. Then it is clear that
\[
\phi(x)=\mathrm{L}_{a}-\mathrm{R}_{a}=\operatorname*{ad}\nolimits_{C}%
a\quad\text{and}\quad\phi(y)=\mathrm{L}_{b}-\mathrm{R}_{b}=\operatorname*{ad}%
\nolimits_{C}b.
\]
Since $a$ and $b$ are Engel elements of $C$, operators $\phi(x)$ and $\phi(y)$
are quasinilpotent elements of $B$. On the other hand, it follows from Lemma
\ref{tens} that $A$ satisfies conditions of Lemma \ref{sum}, whence, by Lemma
\ref{sum}, $\phi(x+y)$ is a quasinilpotent element of $B$. Since
\[
\phi(x+y)=\operatorname*{ad}\nolimits_{C}\left(  a+b\right)  ,
\]
this means that $a+b$ is an Engel element of $C$. But then $a+b$ is an Engel
element of $L$, too.
\end{proof}

For a subspace $N\subset\mathcal{B}(\mathcal{X})$, we will denote by
$\overline{N}$ the closure of $N$ in $\mathcal{B}(\mathcal{X})$.

\begin{corollary}
\label{engelin} Let $L$ be a Lie algebra of compact operators, $L_{0}$ an
$E$-solvable subalgebra of $L$, and let $N$ be the set of all elements of
$L_{0}$ that are Engel elements of $L$. Then $N$ is a closed ideal of $L_{0}$.
\end{corollary}

\begin{proof}
Let $M$ be the set of all elements of $\overline{L_{0}}$ that are Engel
elements of $\overline{L}$. It follows from Theorem \ref{tensor} that $M$ is a
subspace of $\overline{L_{0}}$, whence $M$ is a closed ideal of $\overline
{L_{0}}$ by Lemma \ref{space}. If $N$ is the set of all elements of $L_{0}$
that are Engel elements of $L$, then it is clear that
\[
N=L_{0}\cap M,
\]
and that $L_{0}\cap M$ is a closed ideal of $L_{0}$.
\end{proof}

\begin{remark}
\label{adnil} \emph{Note that }for a finite rank operator\textrm{ }$a$\textrm{
}in a Lie algebra\textrm{ }$L\subset\mathcal{B}\left(  \mathcal{X}\right)
$\textrm{ }the condition of being an Engel element of\textrm{ }$L$\textrm{
\emph{i}}s equivalent to the condition of being an\textrm{ $\operatorname*{ad}%
$\textit{-}}nilpotent element of\textrm{ }$L$\textrm{. \emph{Indeed, since}
}$a$\textrm{ \emph{is an algebraic operator, the operators} $\mathrm{L}$}%
$_{a}$\textrm{ \emph{and} }$\mathrm{R}_{a}$\textrm{ \emph{are also algebraic.
Since they commute, their difference} $\mathrm{L}$}$_{a}-\mathrm{R}_{a}%
$\textrm{ \emph{is algebraic operator whence} }$\operatorname*{ad}a$\textrm{,
\emph{being the restriction of} $\mathrm{L}$}$_{a}-\mathrm{R}_{a}$
\emph{to}\textrm{ }$L$\textrm{,\emph{ is algebraic. But it is evident that}%
}\emph{ a quasinilpotent algebraic operator is nilpotent.}
\end{remark}

\section{\label{_4further}Further Conditions Equivalent to the
Triangularizability of Lie Algebras of Compact Operators}

In this section we will derive some new criteria for the triangularizability
of Lie algebras of compact operators, which, while relevant to the discussion
of graded Lie algebras, are also of independent interest.

For a subset $M$ of a Lie algebra $L$, we define
\[
\mathcal{C}(M)=\{[a,b]\in L:a,b\in M\}.
\]
This set is not additive in general, but is a \textit{Lie multiplicative set},
i.e. it is closed under taking commutators.

It is clear that if $M$ is a Lie multiplicative set in a Lie algebra $L$ then
$\mathcal{C}(M)$ is also Lie multiplicative and $\mathcal{C}(M)\subset M$.
Also, since the map $\operatorname*{ad}:a\mapsto\operatorname*{ad}a$ is a
representation of $L$ on $L$, i.e.
\[
\left[  \operatorname*{ad}a,\operatorname*{ad}b\right]  =\operatorname*{ad}%
\left[  a,b\right]
\]
for every $a,b\in L$, we have that if $M$ is a Lie multiplicative set in $L$,
then $\operatorname*{ad}_{L}M$ is also a Lie multiplicative set in
$\operatorname*{ad}L$.

It is well known \cite[Theorem 2.2.1]{Jac} that \textit{the span of a Lie
multiplicative set of nilpotent operators on a finite dimensional space
consists of nilpotent operators}. The following lemma uses and extends this result.

\begin{lemma}
\label{multi}Let $L$ be a Lie algebra of compact operators, $M$ a Lie
multiplicative set of finite rank operators that are Engel elements of $L$,
and let $L_{0}=\operatorname*{span}M$. Then $L_{0}$ is a subalgebra of $L$
consisting of $\operatorname*{ad}$-nilpotent elements of $L$, and $\left[
L_{0},L_{0}\right]  $ consists of nilpotent finite rank operators. Moreover,
if $M$ consists of nilpotent finite rank operators then so does $L_{0}$.
\end{lemma}

\begin{proof}
It is clear that $L_{0}$ is a subalgebra of $L$. Let $F$ be a finite subset of
$M$, and let $M_{F}$ be a Lie multiplicative set generated by $F$. Then
\[
M_{F}\subset M,
\]
and there is a finite dimensional subspace $\mathcal{Y}$ of $\mathcal{X}$ such
that
\[
\operatorname*{span}M_{F}\subset\{a\in L_{0}:a\mathcal{X}\subset
\mathcal{Y}\}.
\]

Suppose, to the contrary, that $\operatorname*{span}\mathcal{C}(M_{F})$ does
not consist of nilpotent operators. Then there are $a\in\operatorname*{span}%
\mathcal{C}(M_{F})$ and $x\in\mathcal{X}$ such that
\[
ax=x.
\]
It is obvious that $x\in\mathcal{Y}$. Clearly $\mathcal{Y}$ is an invariant
subspace for $\operatorname*{span}M_{F}$, and that $\operatorname*{span}%
\mathcal{C}(M_{F})|\mathcal{Y}$ contains a finite rank operator having
eigenvalue $1$.

Now let $N=\operatorname*{span}M_{F}|\mathcal{Y}$. Then $N$ is a finite
dimensional Lie algebra of operators on $\mathcal{Y}$, because $\mathcal{Y}$
is finite dimensional. Since $M_{F}$ consists of Engel elements of $L$,
$\operatorname*{ad}_{N}M_{F}|\mathcal{Y}$ is a Lie multiplicative set of
nilpotent operators on the finite dimensional space $N$, whence so is
$\operatorname*{span}\operatorname*{ad}_{N}M_{F}|\mathcal{Y}$ by \cite[Theorem
2.2.1]{Jac}. Since
\[
\operatorname*{span}\operatorname*{ad}\nolimits_{N}M_{F}|\mathcal{Y}%
=\operatorname*{ad}\nolimits_{N}\operatorname*{span}M_{F}|\mathcal{Y},
\]
it follows that the Lie algebra $\operatorname*{span}M_{F}|\mathcal{Y}$ is
Engel, and hence is triangularizable by Theorem \ref{redcri}. In particular,
$\operatorname*{span}\mathcal{C}(M_{F})|\mathcal{Y}$ consists of nilpotent
operators, which gives a contradiction.

Thus $\operatorname*{span}\mathcal{C}(M_{F})$ consists of nilpotent operators.
Since $F$ is arbitrary, we obtain in fact that $\operatorname*{span}%
\mathcal{C}(M)$ consists of nilpotent finite rank operators. Since
\[
\operatorname*{span}\mathcal{C}(M)=\left[  L_{0},L_{0}\right]  ,
\]
$L_{0}$ is $E$-solvable and triangularizable by Theorem \ref{redcri}.

By Theorem \ref{tensor}, $\operatorname*{span}M$ consists of Engel elements of
$L$, but $\operatorname*{span}M=L_{0}$. By Remark \ref{adnil}, $L_{0}$
consists of $\operatorname*{ad}$-nilpotent elements of $L$. This completes the
proof of the basic part of the lemma.

If $M$ consists of nilpotent finite rank operators then, applying the
subadditivity of the spectral radius on $L_{0}$ (see Theorem \ref{redcri}), we
obtain that $L_{0}$ is Volterra, and hence it consists of nilpotent finite
rank operators.
\end{proof}

We need the following technical lemmas.

\begin{lemma}
\label{cl71} Let $\omega:L\longrightarrow M$ be a bounded homomorphism with
dense image from a normed Lie algebra $L$ into a Lie subalgebra $M$ of a
compact Banach algebra $B$. Then
\[
\rho(\left(  \operatorname*{ad}\omega(a)\right)  \left(  \operatorname*{ad}%
\omega\left(  b\right)  \right)  )\leq\rho(\left(  \operatorname*{ad}a\right)
\left(  \operatorname*{ad}b\right)  )
\]
for every $a,b\in L$.
\end{lemma}

\begin{proof}
One can assume that $M$ is closed in $B$. Let $T=\left(  \operatorname*{ad}%
\omega\left(  a\right)  \right)  \left(  \operatorname*{ad}\omega\left(
b\right)  \right)  $ and $S=\left(  \mathrm{L}_{\omega\left(  a\right)
}-\mathrm{R}_{\omega\left(  a\right)  }\right)  \left(  \mathrm{L}%
_{\omega\left(  b\right)  }-\mathrm{R}_{\omega\left(  b\right)  }\right)  $.
Note that
\[
S=S_{1}-S_{2}=\left(  \mathrm{L}_{\omega\left(  a\right)  \omega\left(
b\right)  }+\mathrm{R}_{\omega\left(  b\right)  \omega\left(  a\right)
}\right)  -\left(  \mathrm{L}_{\omega\left(  a\right)  }\mathrm{R}%
_{\omega\left(  b\right)  }+\mathrm{L}_{\omega\left(  b\right)  }%
\mathrm{R}_{\omega\left(  a\right)  }\right)  ,
\]
where the operator
\[
S_{1}=\mathrm{L}_{\omega\left(  a\right)  \omega\left(  b\right)  }%
+\mathrm{R}_{\omega\left(  b\right)  \omega\left(  a\right)  }%
\]
has countable spectrum on $B$ since $\omega\left(  a\right)  \omega\left(
b\right)  $ and $\omega\left(  b\right)  \omega\left(  a\right)  $ have
countable spectra as compact elements of $B$, and $\sigma(S_{1})\subset
\sigma\left(  \omega\left(  a\right)  \omega\left(  b\right)  \right)
+\sigma\left(  \omega\left(  b\right)  \omega\left(  a\right)  \right)  $, and
the operator
\begin{align*}
S_{2}  &  =\mathrm{L}_{\omega\left(  a\right)  }\mathrm{R}_{\omega\left(
b\right)  }+\mathrm{L}_{\omega\left(  b\right)  }\mathrm{R}_{\omega\left(
a\right)  }\\
&  =\mathrm{L}_{\omega\left(  a\right)  +\omega\left(  b\right)  }%
\mathrm{R}_{\omega\left(  a\right)  +\omega\left(  b\right)  }-\mathrm{L}%
_{\omega\left(  a\right)  }\mathrm{R}_{\omega\left(  a\right)  }%
-\mathrm{L}_{\omega\left(  b\right)  }\mathrm{R}_{\omega\left(  b\right)  }%
\end{align*}
is compact on $B$ by definition of compact Banach algebras. Therefore $S$ has
countable spectrum as a compact perturbation of an operator with countable
spectrum (for instance see \cite[Proposition 3.23]{ShT2005}). Then $T$ has
countable spectrum as the restriction of $S$ to $M$.

For every $c\in L$, we have that%
\begin{align*}
\left\Vert \left(  \left(  \operatorname*{ad}\omega\left(  a\right)  \right)
\left(  \operatorname*{ad}\omega\left(  b\right)  \right)  \right)  ^{n}%
\omega\left(  c\right)  \right\Vert ^{1/n} &  =\left\Vert \omega\left(
\left(  \left(  \operatorname*{ad}a\right)  \left(  \operatorname*{ad}%
b\right)  \right)  ^{n}c\right)  \right\Vert ^{1/n}\\
&  \leq\left(  \left\Vert \omega\right\Vert \left\Vert c\right\Vert \right)
^{1/n}\left\Vert \left(  \left(  \operatorname*{ad}a\right)  \left(
\operatorname*{ad}b\right)  \right)  ^{n}\right\Vert ^{1/n}%
\end{align*}
for each integer $n$, whence
\[
\lim\sup\left\Vert T^{n}x\right\Vert ^{1/n}\leq\rho(\left(  \operatorname*{ad}%
a\right)  \left(  \operatorname*{ad}b\right)  )
\]
for every $x$ in the image of $\omega$. Recall, as $T$ has countable spectrum,
that the subspace $\{x\in M:\lim\sup\left\Vert T^{n}x\right\Vert ^{1/n}\leq
r\}$ is in fact closed (see, for instance, \cite[Proposition 3.3]{ShT2005})
for every $r\geq0$. Since one can consider $M$ as a Banach space and
$\omega\left(  L\right)  $ is dense in $M$, we obtain by (\ref{f_muller})
that
\[
\rho(T)=\sup_{x\in M}\lim\sup\left\Vert T^{n}x\right\Vert ^{1/n}\leq
\rho(\left(  \operatorname*{ad}a\right)  \left(  \operatorname*{ad}b\right)
).
\]

\end{proof}

\begin{remark}
\label{radius} $\emph{The}$ $\emph{same}$ $\emph{argument}$ $\emph{as}$
$\emph{in}$ $\emph{Lemma}$ $\mathrm{\ref{cl71}}$ $\emph{proves}$ $\emph{that}$
if $\omega:L\longrightarrow M$ is a bounded homomorphism with dense image from
a normed Lie algebra $L$ into a Banach Lie algebra $M$ then%
\[
\rho(\operatorname*{ad}\omega(a))\leq\rho(\operatorname*{ad}a)
\]
for every $a\in L$ such that $\operatorname*{ad}\omega(a)$ is an operator with
countable spectrum.

$\emph{This}$ $\emph{shows}$ $\emph{for}$ $\emph{instance}$ $\emph{that}$ the
property of being Engel is preserved for compact operators induced on a
quotient $\mathcal{V}$ of invariant closed subspaces; $\emph{for}$
$\emph{this,}$ $\emph{it}$ $\emph{suffices}$ $\emph{to}$ $\emph{take}$
$L\subset\mathcal{K}(\mathcal{X})$, $M=\overline{L|\mathcal{V}}$ $\emph{and}$
$\omega(a)=a|\mathcal{V}$ $\emph{for}$ $\emph{every}$ $a\in L$.
\end{remark}

\begin{lemma}
\label{cl72} Let $L$ be a Lie algebra of compact operators on $\mathcal{X}$,
and let $y\in L\cap\mathcal{F}(\mathcal{X})$ be such that $\rho(\left(
\operatorname*{ad}y\right)  \left(  \operatorname*{ad}b\right)  )=0$ for every
$b\in L$. Then $\{y,z\}$ is triangularizable for every $z\in L\cap
\mathcal{F}(\mathcal{X})$.
\end{lemma}

\begin{proof}
If not then there are closed subspaces $Z\subset Y$ invariant for both $y$ and
$z$, such that the restriction of $\{y,z\}$ to the quotient $\mathcal{V}=Y/Z$
is irreducible and $\dim V>1$. It is clear that $y|\mathcal{V}$ and
$z|\mathcal{V}$ are non-scalar.

Let $N$ be the Lie algebra generated by $y$ and $z$. Since $y|\mathcal{V}$ and
$z|\mathcal{V}$ are of finite rank, it is clear that $\mathcal{V}$ is finite
dimensional. It follows from Lemma \ref{cl71} that
\[
\rho(\left(  \operatorname*{ad}y|\mathcal{V}\right)  \left(
\operatorname*{ad}b|\mathcal{V}\right)  )=0
\]
for every $b\in N$, whence in particular we have that
\[
\left\langle y|\mathcal{V},N|\mathcal{V}\right\rangle =0,
\]
where the Killing form is taken relative to the Lie algebra $N|\mathcal{V}$
which is finite dimensional. Let
\[
I=\left\{  a\in N|\mathcal{V}:\left\langle a,N|\mathcal{V}\right\rangle
=0\right\}  .
\]
Since
\[
\left\langle \left[  a,b\right]  ,N|\mathcal{V}\right\rangle =\left\langle
a,\left[  b,N|\mathcal{V}\right]  \right\rangle =0
\]
for every $a\in I$ and $b\in N|\mathcal{V}$, the subspace $I$ is in fact an
ideal of $N|\mathcal{V}$. Since also%
\[
\left\langle I,\left[  I,I\right]  \right\rangle =0,
\]
$I$ is solvable by Lemma \ref{cartan}. Since $I$ contains a non-scalar
operator, $N|\mathcal{V}$ is reducible by Theorem \ref{solid}, which is a contradiction.
\end{proof}

\begin{lemma}
\label{claim3} Let $L$ be a Lie algebra of compact operators, $B=\mathcal{A}%
(L)/\operatorname*{rad}\mathcal{A}(L)$, $\theta:\mathcal{A}(L)\longrightarrow
B$ a canonical map, and let $M$ be the closure of $\theta(L)$ in $B$. If $M$
is Engel, then $L$ is $E$-solvable.
\end{lemma}

\begin{proof}
Let us show that $B$ is commutative. Since $B$ is semisimple, we have only to
show that each strictly irreducible representation $\tau$ of $B$ acts on
one-dimensional space.

Since
\[
\tau(\mathcal{E}_{0}(\operatorname*{ad}z))\subset\mathcal{E}_{0}%
(\operatorname*{ad}\tau(z))
\]
for every $z\in M$, we obtain from our hypothesis that $\operatorname*{ad}%
\tau(z)$ is locally quasinilpotent on $\tau(M)$, and hence that it is
quasinilpotent. By \cite[Theorem 5.1]{Ale}, $\tau(M)$ is an Engel Lie algebra
of compact operators, whence $\tau(B)$ is triangularizable by Theorem
\ref{redcri}. Since $\tau(B)$ is irreducible, $\tau$ is one-dimensional.

Thus $L$ is triangularizable by \cite{Mu2}, and is $E$-solvable by Theorem
\ref{redcri}.
\end{proof}

Recall \cite[Theorems 5.19 and 5.5]{ShT2005} that there exists a largest Engel
(respectively, Volterra) ideal in any Lie algebra of compact operators.

\begin{lemma}
\label{volterraengel}Let $L$ be a closed Lie algebra of compact operators on
$\mathcal{X}$, and let $E$ be the largest Engel ideal of $L$. If all nilpotent
finite rank operators in $L$ belong to $E$, then $L/E$ is Engel, $L$ is
$E$-solvable and therefore triangularizable.
\end{lemma}

\begin{proof}
For a bounded operator $x$ on a Banach space, let $p_{\lambda,r}(x)$ be the
Riesz projection of $x$ corresponding to the part of $\sigma(x)$ enclosed by
the circle $\Omega_{\lambda,r}$ of radius $r$ centered at $\lambda
\in\mathbb{C}$ whenever this circle lies into the resolvent set of $x$. Recall
that
\begin{equation}
p_{\lambda,r}(x)=\frac{1}{2\pi i}\int_{\Omega_{\lambda,r}}(\mu-x)^{-1}d\mu.
\label{projection}%
\end{equation}
In the following claims we will deal with operators (of adjoint
representation) with countable spectra. Let $\varphi:L\longrightarrow L/E$ be
\textbf{the} canonical map.

\begin{claim}
$L/E$ is an Engel Lie algebra.
\end{claim}

If not, there are an element $a\in L$ and a non-zero $\lambda\in\mathbb{C}$
such that $\lambda$ is an isolated point of $\sigma(\operatorname*{ad}%
\varphi(a))$. Note that
\[
\operatorname*{ad}\varphi(a)=\left(  \operatorname*{ad}a\right)  |\varphi(L).
\]
It is clear that there are $p_{\lambda,r}(\operatorname*{ad}a)$ and
$p_{\lambda,r}(\operatorname*{ad}\varphi(a))$, where the circle $\Omega
_{\lambda,r}$ encloses only $\{\lambda\}$ in $\sigma(\operatorname*{ad}%
\varphi(a))$ and $r<\left\vert \lambda\right\vert $. Since $\lambda\in$
$\sigma(\operatorname*{ad}\varphi(a))$, we have that
\[
p_{\lambda,r}(\operatorname*{ad}\varphi(a))\neq0.
\]
By using the fact that the resolvent set of an operator with countable
spectrum is connected, we obtain that
\begin{equation}
(\mu-\operatorname*{ad}a)^{-1}|\varphi(L)=(\mu-\left(  \operatorname*{ad}%
a\right)  |\varphi(L))^{-1} \label{inverse}%
\end{equation}
for every $\mu$ in the intersection of resolvent sets of $a$ and
$a|\varphi(L)$. It follows from (\ref{projection}) (under
$x=\operatorname*{ad}a$) and (\ref{inverse}) that
\[
p_{\lambda,r}(\operatorname*{ad}\varphi(a))=p_{\lambda,r}\left(  \left(
\operatorname*{ad}a\right)  |\varphi(L)\right)  =p_{\lambda,r}\left(
\operatorname*{ad}a\right)  |\varphi(L)
\]
Note that
\[
p_{\lambda,r}\left(  \operatorname*{ad}a\right)  L=\mathcal{E}_{\lambda
,r}(\operatorname*{ad}a):=\left\{  x\in L:\lim\sup\left\Vert \left(
\operatorname*{ad}a-\lambda\right)  ^{n}x\right\Vert ^{1/n}\leq r\right\}  .
\]
As $r<\left\vert \lambda\right\vert $, it follows from \cite[Corollary
3.12]{ShT2005} that $\mathcal{E}_{\lambda,r}(\operatorname*{ad}a)$ consists of
nilpotent finite rank operators. Then we have that
\[
\mathcal{E}_{\lambda,r}(\operatorname*{ad}a)\subset E.
\]
This shows that $p_{\lambda,r}\left(  \operatorname*{ad}a\right)  L\subset E$
and
\[
p_{\lambda,r}(\operatorname*{ad}\varphi(a))=p_{\lambda,r}\left(
\operatorname*{ad}a\right)  |\varphi(L)=0,
\]
a contradiction.

\begin{claim}
$L$ is $E$-solvable and triangularizable.
\end{claim}

Let $I$ be a proper closed ideal of $L$, and let $\omega:L\longrightarrow L/I$
be a canonical map. It is clear that $\operatorname*{ad}\omega\left(
a\right)  $ is an operator with countable spectrum for every $a\in L$. If $I$
doesn't contain $E$ then it follows from Remark \ref{radius} that
$\omega\left(  E\right)  $ is a non-zero Engel ideal of $L/I$. Its closure is
also Engel in virtue of continuity of the spectral radius on operators with
countable spectra.

Suppose now that $E\subset I$. Then there is a natural bounded epimorphism
$L/E\longrightarrow L/I,$ whence $L/I$ is clearly Engel.

Therefore $L$ is $E$-solvable. Then $L$ is triangularizable by Theorem
\ref{redcri}.
\end{proof}

Now we are in a position to obtain the following result which complements
Theorem \ref{redcri}.

\begin{theorem}
\label{crit} The following conditions can be added to the list of conditions
equivalent to triangularizability in Theorem $\ref{redcri}$, for a Lie algebra
$L$ of compact operators.

\begin{enumerate}
\item[(7)] The set $\mathcal{C}(L)$ of all commutators of $L$ consists of
Engel elements of $L$.

\item[(8)] $\mathcal{C}(L)$ is Volterra.

\item[(9)] $\rho(\left(  \operatorname*{ad}a\right)  \left(
\operatorname*{ad}b\right)  )\leq\rho(\operatorname*{ad}a)\rho
(\operatorname*{ad}b)$ for every $a,b\in L$.

\item[(10)] $\rho(ab)\leq\rho(a)\rho(b)$ for every $a,b\in L$.

\item[(11)] $a+b$ is an Engel element of $\overline{L}$ for every
$\operatorname*{ad}$-nilpotent elements $a,b$ of $\overline{L}$.

\item[(12)] $a+b$ is nilpotent for every nilpotent finite rank operators $a,b$
in $\overline{L}$.

\item[(13)] $[a,b]$ is an Engel element of $\overline{L}$ for every
$\operatorname*{ad}$-nilpotent elements $a,b$ of $\overline{L}$.

\item[(14)] The set of all nilpotent finite rank operators in $\overline{L}$
is Lie multiplicative.
\end{enumerate}
\end{theorem}

\begin{proof}
Let $B=\mathcal{A}(L)/\operatorname*{rad}\mathcal{A}(L)$, $\theta
:\mathcal{A}(L)\longrightarrow B$ a canonical map, and let $M$ be the closure
of $\theta(L)$ in $B$.

Recall that $(2)$ of Theorem $\ref{redcri}$ means that $L$\textit{ is
triangularizable}, and we refer to it simply as $(2)$. We start with an easy
part of the implications.

\begin{itemize}
\item $(2)$ implies either of $(7)$, $(8)$, $(9)$, $(10)$, $(11)$, $(12)$,
$(13)$, and $(14)$.
\end{itemize}

Indeed, the triangularizability of $L$ implies (9) by \cite[Theorem 5.12 and
Corrolary 5.13]{ShT2005}, (11) by \cite[Corollary 5.13 and Theorem
5.19]{ShT2005}, and (13) by \cite[Theorem 5.19]{ShT2005}. The other
implications are obvious, because if $L$ is triangularizable, then
$\mathcal{A}(L)$ is commutative modulo the Jacobson radical, and hence has the
spectral properties required for these implications.

Now we must show the reverse implications.

\begin{itemize}
\item $(7)\Longrightarrow(2)$.
\end{itemize}

Let
\[
I=\mathcal{C}(L)\cap\mathcal{F}(\mathcal{X}).
\]
Then $I$ is a Lie multiplicative set in $L$. By Lemma \ref{cl71}, the property
being Engel is preserved for compact operators induced on a quotient of
invariant subspaces, whence for the implication it is sufficient only to show
that $L$ is reducible.

Suppose, to the contrary, that $L$ is irreducible. Then $L$ is not Engel, and
hence there is a non-zero nilpotent finite rank operator $a$ of $L$ by Theorem
\ref{engel}. If $[a,L]=0$, then $L_{0}$ is reducible by Lomonosov's Theorem
\cite{Lom}, which gives a contradiction. Therefore $[a,L]\neq0$ , whence%
\[
I\neq0,
\]
and moreover, $I$ is non-scalar. Since, by hypothesis, $\mathcal{C}(L)$
consists of Engel elements of $L$ so does the set $I$. By Lemma \ref{multi},
$\operatorname*{span}I$ consists of Engel elements of $L$. In particular,
$\operatorname*{span}I$ is Engel. But $\operatorname*{span}I$ is an ideal of
$L$, whence $L$ is reducible by Theorem \ref{solid}, a contradiction.

\begin{itemize}
\item $(8)\Longrightarrow(2).$
\end{itemize}

It is clear that $(8)\Longrightarrow(7)$ is true, whence $(8)\Longrightarrow
(2)$ holds in virtue of $(7)\Longrightarrow(2)$.

\begin{itemize}
\item $(9)\Longrightarrow(2).$
\end{itemize}

Assume that $(9)$ holds. It follows from Lemma \ref{claim3} that it is
sufficient to show that $M$ is Engel.

\begin{claim}
\label{cl73} Let $y$ be a quasinilpotent element in $M$. Then $\rho(\left(
\operatorname*{ad}y\right)  \left(  \operatorname*{ad}z\right)  )=0$ for every
$z\in M$.
\end{claim}

Suppose first that $z=\theta(b)$ for some $b\in L$. Let $(a_{k})$ be a
sequence of elements of $L$ such that
\[
\theta(a_{k})\rightarrow y
\]
as $k\rightarrow\infty$. Then we have that
\[
\rho\left(  a_{k}\right)  =\rho\left(  \theta(a_{k})\right)  \rightarrow0,
\]
and hence that%
\[
\rho\left(  \operatorname*{ad}a_{k}\right)  \leq2\rho\left(  a_{k}\right)
\rightarrow0,
\]
as $k\rightarrow\infty$. It follows from Lemma \ref{cl71} and (9) that
\[
\rho(\left(  \operatorname*{ad}\theta(a_{k})\right)  \left(
\operatorname*{ad}\theta(b)\right)  )\leq\rho(\left(  \operatorname*{ad}%
a_{k}\right)  \left(  \operatorname*{ad}b\right)  )\leq\rho\left(
\operatorname*{ad}a_{k}\right)  \rho\left(  \operatorname*{ad}b\right)
\rightarrow0,
\]
and since $\left(  \operatorname*{ad}\theta(a_{k})\right)  \left(
\operatorname*{ad}\theta(b)\right)  \rightarrow\left(  \operatorname*{ad}%
y\right)  \left(  \operatorname*{ad}z\right)  $ as $k\rightarrow\infty$, we
obtain that
\[
\rho(\left(  \operatorname*{ad}y\right)  \left(  \operatorname*{ad}z\right)
)=0
\]
for every $z$ in $\theta(L)$. By continuity of $\rho$ on elements with
countable spectra, this equality holds for every $z\in M$.

\begin{claim}
\label{cl74} Let $y\in M\cap\operatorname*{soc}B$ be nilpotent. Then $\left[
y,z\right]  $ is a nilpotent element of $\operatorname*{soc}B$ for every $z\in
M\cap\operatorname*{soc}B$.
\end{claim}

Since the spectrum of $\left[  y,z\right]  $ equals the union of
$\sigma\left(  \tau\left(  \left[  y,z\right]  \right)  \right)  $, when
$\tau$ runs over strictly irreducible representations of $B$, it is sufficient
to show that $\sigma\left(  \left[  \tau\left(  y\right)  ,\tau\left(
z\right)  \right]  \right)  =0$ for an arbitrary $\tau$.

Since $\tau(B)$ consists of compact operators by \cite[Theorem 5.1]{Ale},
$\tau(M)$ is a Lie algebra of compact operators, and also $\tau\left(
y\right)  $ and $\tau\left(  z\right)  $ are finite rank operators. It follows
by Claim \ref{cl73} and Lemma \ref{cl71} that%
\[
\rho(\left(  \operatorname*{ad}\tau\left(  y\right)  \right)  \left(
\operatorname*{ad}\tau\left(  b\right)  \right)  )=0
\]
for every $b\in M$. By Lemma \ref{cl72}, $\left\{  \tau\left(  y\right)
,\tau\left(  z\right)  \right\}  $ is triangularizable, whence
\[
\rho\left(  \left[  \tau\left(  y\right)  ,\tau\left(  z\right)  \right]
\right)  =0.
\]
This proves the claim.

\begin{claim}
\label{cl75} Let $N$ be the set of all nilpotent elements in $M\cap
\operatorname*{soc}B$. Then $N$ is an ideal of $M\cap\operatorname*{soc}B$.
\end{claim}

It is sufficient to show that $\tau(\operatorname*{span}N)$ consists of
nilpotent operators for an arbitrary strictly irreducible representation
$\tau$ of $B$.

Since $N$ is Lie multiplicative in $M$ by Claim \ref{cl74}, $\tau(N)$ is a Lie
multiplicative set of nilpotent finite rank operators in $\tau(M)$, whence so
is $\operatorname*{span}\tau(N)$ by Lemma \ref{multi}. This proves the claim.

\begin{claim}
\label{cl76} $N=0$.
\end{claim}

It is sufficient to show that $\tau(N)=0$ for an arbitrary strictly
irreducible representation $\tau$ of $B$.

Suppose, to the contrary, that $\tau(N)\neq0$ for some $\tau$. Since $\tau(N)$
is a Volterra ideal of the Lie algebra $\tau\left(  M\cap\operatorname*{soc}%
B\right)  $,
\[
\tau(N)\subset R,
\]
where $R$ is the largest $E$-solvable ideal of $\tau\left(  M\cap
\operatorname*{soc}B\right)  $. Since $R$ is inner-characteristic by Theorem
\ref{inner} and $\tau\left(  M\cap\operatorname*{soc}B\right)  $ is an ideal
of the Lie algebra $\tau\left(  M\right)  $, $R$ is an $E$-solvable ideal of
$\tau\left(  M\right)  $, whence $\tau\left(  M\right)  $ is reducible by
Theorem \ref{solid}, which is a contradiction.

\begin{claim}
\label{cl77} $M$ is an Engel Lie algebra.
\end{claim}

Indeed, if not, there are non-zero $x,y\in M$, and a non-zero isolated point
$\lambda$ in the spectrum $\sigma(\operatorname*{ad}x)$, such that
\[
y\in\mathcal{E}_{\lambda}(\operatorname*{ad}x).
\]
It follows that $y$ is a nilpotent element of $M\cap\operatorname*{soc}B$ by
Lemma \ref{socle}, whence
\[
y\in N.
\]
But $N=0$ by Claim \ref{cl76}, a contradiction. This proves the claim.

Now $L$ is triangularizable by Lemma \ref{claim3}.

\begin{itemize}
\item $(10)\Longrightarrow(2)$.
\end{itemize}

Assume that (10) holds.

First of all, note that $(10)$\textit{ holds for every }$a,b\in M$.

Indeed, since $\rho(\theta(a))=\rho(a)$ for every $a\in L$, and since the
spectral radius is continuous on elements with countable spectra, (10) holds
for elements of $M$.

Further, we claim that $M$ \textit{is an Engel Lie algebra}.

If not, there are non-zero $x,y\in M$ such that
\[
y\in\mathcal{E}_{1}(\operatorname*{ad}x).
\]
Since $y$ is nilpotent by Lemma \ref{socle}, we have by (10) (reformulated for
elements of $M$) that%
\[
\rho(yz)=0
\]
for every $z\in M$.

Let $\tau$ be a strictly irreducible representation of $B$. By \cite[Theorem
5.1]{Ale}, $\tau(B)$ consists of compact operators. Since $\tau(y)\in
\mathcal{E}_{1}(\operatorname*{ad}\tau(x)),$we have that $\tau(y)$ is a finite
rank operator by \cite[Theorem 3]{Woj}, and since
\[
\rho(\tau(z))\leq\rho(z)
\]
for every $z\in B$, we have that $\tau(yM)$ consists of nilpotent finite rank
operators, whence
\[
\operatorname*{tr}(\tau(y)\tau(M))=0.
\]
If $\tau(y)\neq0$, it follows from Lemma \ref{carfin} that $\tau(M)$ is
reducible, which would be a contradiction.

Then $\tau(y)=0$ what implies, since $\tau$ is arbitrary, that $y=0$, a contradiction.

This proves that $M$ is Engel. Now $L$ is triangularizable by Lemma
\ref{claim3}.

\begin{itemize}
\item $(11)\Longrightarrow(2).$
\end{itemize}

Let $E$ be the set of all Engel elements of $\overline{L}$. By Remark
\ref{adnil}, $E\cap\mathcal{F}(\mathcal{X})$ consists of $\operatorname*{ad}%
$-nilpotent elements of $\overline{L}$, whence $E\cap\mathcal{F}(\mathcal{X})$
is additive by hypothesis. By Lemma \ref{space}, $E\cap\mathcal{F}%
(\mathcal{X})$ is an ideal of $L$, whence $E\cap\mathcal{F}(\mathcal{X})$ is
an Engel ideal of $L$. Then the set of all nilpotent finite rank operators in
$\overline{L}$ is contained in the largest Engel ideal of $\overline{L}$. By
Lemma \ref{volterraengel}, $\overline{L}$ is triangularizable.

\begin{itemize}
\item $(12)\Longrightarrow(2).$
\end{itemize}

Let $F$ be the set of all nilpotent finite rank operators in $\overline{L}$.
Since $F$ is a subspace of $L$, we have by Lemma \ref{space} that $F$ is an
ideal of $L$, whence it is contained in the largest Engel ideal of
$\overline{L}$. By Lemma \ref{volterraengel}, $\overline{L}$ is triangularizable.

\begin{itemize}
\item $(13)\Longrightarrow(2).$
\end{itemize}

Let $E$ be the set of all Engel elements of $\overline{L}$. By Remark
\ref{adnil}, $E\cap\mathcal{F}(\mathcal{X})$ consists of $\operatorname*{ad}%
$-nilpotent elements of $\overline{L}$, and is therefore Lie multiplicative by
hypothesis. By Lemma \ref{multi},
\[
\operatorname*{span}E\cap\mathcal{F}(\mathcal{X})\subset E\cap\mathcal{F}%
(\mathcal{X}),
\]
whence $E\cap\mathcal{F}(\mathcal{X})$ is an ideal of $\overline{L}$ by Lemma
\ref{space}. Then $E\cap\mathcal{F}(\mathcal{X})$ is an Engel ideal of
$\overline{L}$, and the set of all nilpotent finite rank operators in
$\overline{L}$ is contained in the largest Engel ideal of $\overline{L}$. By
Lemma \ref{volterraengel}, $\overline{L}$ is triangularizable.

\begin{itemize}
\item $(14)\Longrightarrow(2).$
\end{itemize}

Let $F$ be the set of all nilpotent finite rank operators in $\overline{L}$.
By hypothesis, $F$ is Lie multiplicative. By Lemma \ref{multi},
\[
\operatorname*{span}F\subset F.
\]
By Lemma \ref{space}, $F$ is an ideal of $L$, whence it is contained in the
largest Engel ideal of $\overline{L}$. By Lemma \ref{volterraengel},
$\overline{L}$ is triangularizable.
\end{proof}

\begin{remark}
\label{addcrit} \emph{Using the continuity of the spectral radius on operators
with countable spectra, it is easy to see that} $(5)$ \emph{and }$(6)$
\emph{of Theorem} $\ref{redcri}$ \emph{imply} $(11)$ \emph{and} $(12)$
\emph{of Theorem} $\ref{crit}$, \emph{respectively. So, for applications to
triangularizability, conditions} $(11)$ \emph{and} $(12)$ \emph{are more
effective than} $(5)$\emph{ and }$(6)$. \emph{Similarly, the following
conditions}

\begin{enumerate}
\item[(15)] The spectral radius is Lie submultiplicative on
$\operatorname*{ad}L.$

\item[(16)] The spectral radius is Lie submultiplicative on $L$ (i.e.
$\rho\left(  \left[  a,b\right]  \right)  \leq2\rho\left(  a\right)
\rho\left(  b\right)  $ for every $a,b\in L$).
\end{enumerate}

\noindent\emph{are stronger than} $(13)$ \emph{and }$(14)$,
\emph{respectively, and clearly follow from} $(2)$ \emph{of Theorem}
$\ref{redcri}$; \emph{as a result they are equivalent to the
triangularizability of a Lie algebra} $L$ \emph{of compact operators. Clearly}
$(13)$ \emph{and} $(14)$ \emph{are more effective tools for establishing
triangularizability than} $(15)$ \emph{and }$(16)$.
\end{remark}

Note that the conditions in Theorems \ref{redcri} and \ref{crit}, and in
Remark \ref{addcrit}, are arranged in pairs; odd conditions are expressed in
terms of normed Lie algebras, while even conditions are expressed in terms of
operator algebras.

\part{\label{__2graded}Subgraded Lie Algebras of Compact Operators}

Our aim in this part is to determine which of quasinilpotence conditions on
homogeneous elements of a subgraded Lie algebra of compact operators are
equivalent to triangularizability. The main result is Theorem \ref{main1}.

\section{\label{_5tensor} Graded Ampliations of Subgraded Lie Algebras}

Let $\Gamma$ be a finite commutative group, and let $\pi$ be a faithful
representation of $\Gamma$ on a finite dimensional space $\mathcal{Y}$ (i.e.
an injective homomorphism of $\Gamma$ into the group of invertible operators
on $\mathcal{Y}$), and let $L$ be a $\Gamma$-subgraded Lie algebra of
operators in $\mathcal{K}(\mathcal{X})$. Instead of working directly with $L$,
it is convenient to consider a related Lie algebra $M:=$ $L^{\pi}$ of
operators in $\mathcal{K}(\mathcal{Z})$ where $\mathcal{Z}=\mathcal{X}%
\otimes\mathcal{Y}$, and $M$ is the linear space generated by subspaces
\[
M_{\gamma}=L_{\gamma}\otimes\pi(\gamma),
\]
$\gamma\in\Gamma$.

It is clear that $M=\sum_{\gamma\in\Gamma}M_{\gamma}$ where the sum is direct,
and that%
\[
\lbrack M_{\alpha},M_{\beta}]\subset M_{\alpha+\beta}%
\]
for every $\alpha,\beta\in\Gamma$, so $M$ is a $\Gamma$-graded Lie algebra.

Setting
\[
f_{\pi}(a\otimes\pi(\gamma))=a
\]
for all $\gamma\in\Gamma$, $a\in L_{\gamma}$ and extending it linearly to $M$,
we obtain a continuous homomorphism $f_{\pi}$ of $M$ onto $L$. Clearly
$f_{\pi}$ \textit{is injective if and only if }$L$\textit{ is graded.} We will
refer to $L^{\pi}$ constructed in this way as a \textit{graded ampliation} of
$L$.

For a $\mathbb{Z}_{n}$-subgraded Lie algebra $L$, it is not difficult to
realize $L^{\pi}$ by concrete $n\times n$-matrices with entries in $L$. For
$n=2$, this is shown in the following example.

\begin{example}
For a $\mathbb{Z}_{2}$-subgraded Lie algebra $L=L_{0}+L_{1}$, define $\pi$ by
\[
\pi(0)=%
\begin{pmatrix}
1 & 0\\
0 & 1
\end{pmatrix}
,\text{\quad}\pi(1)=%
\begin{pmatrix}
0 & 1\\
1 & 0
\end{pmatrix}
.
\]
Then $L^{\pi}$ is equal to $M=M_{0}+M_{1}$ up to an isomorphism, where%
\[
M_{0}=\left\{
\begin{pmatrix}
a & 0\\
0 & a
\end{pmatrix}
:a\in L_{0}\right\}  ,\quad M_{1}=\left\{
\begin{pmatrix}
0 & b\\
b & 0
\end{pmatrix}
:b\in L_{1}\right\}  .
\]
In this case,%
\[
f_{\pi}\left(
\begin{pmatrix}
a & b\\
b & a
\end{pmatrix}
\right)  =a+b.
\]

\end{example}

In what follows, $\Gamma$ is assumed to be finite, unless specifically
indicated otherwise.

\begin{lemma}
\label{maptri}Let $L$ be a $\Gamma$-subgraded Lie algebra of compact
operators. If $L^{\pi}$ is Engel (or $E$-solvable) then so is $L$.
\end{lemma}

\begin{proof}
For an arbitrary homomorphism of Lie algebras, we have
\[
(\operatorname*{ad}f_{\pi}{(a)})f_{\pi}=f_{\pi}\operatorname*{ad}a
\]
for every $a\in L^{\pi}$, so it follows that
\[
(\operatorname*{ad}f_{\pi}{(a)})^{n}f_{\pi}=f_{\pi}(\operatorname*{ad}a)^{n}%
\]
for all $n$.

If $L^{\pi}$ is Engel then $\Vert(\operatorname*{ad}a)^{n}\Vert^{1/n}%
\rightarrow0$ as $n\rightarrow\infty$, for every $a\in L^{\pi}$, whence
\[
\Vert(\operatorname*{ad}f_{\pi}{(a)})^{n}f_{\pi}\Vert^{1/n}\rightarrow0.
\]
Since $f_{\pi}(L^{\pi})=L$, this shows that
\[
\Vert(\operatorname*{ad}x)^{n}y\Vert^{1/n}\rightarrow0
\]
as $n\rightarrow\infty$, for every $x,y\in L$. As $\operatorname*{ad}x$ has
countable spectrum, $\operatorname*{ad}x$ is decomposable. So the local
quasinilpotence of $\operatorname*{ad}x$ on $L$ implies that
$\operatorname*{ad}x$ is quasinilpotent, and hence $L$ is Engel.

Now if $L^{\pi}$ is $E$-solvable then $[L^{\pi},L^{\pi}]$ is Volterra, whence
$[L^{\pi},L^{\pi}]$ is Engel. It follows from above that $f_{\pi}([L^{\pi
},L^{\pi}])$ is Engel. Since
\[
\lbrack L,L]=f_{\pi}([L^{\pi},L^{\pi}]),
\]
$L$ is $E$-solvable by Theorem \ref{redcri}.
\end{proof}

\begin{theorem}
\label{cons} Let $L$ be a $\Gamma$-subgraded Lie algebra of compact operators.
If the union of the closures of all components $L_{\gamma}$ doesn't contain
non-zero finite rank operators, then $L$ is Engel.
\end{theorem}

\begin{proof}
Replacing $L$ by $\sum_{\gamma\in\Gamma}\overline{L_{\gamma}}$ if necessary,
we may assume that each $L_{\gamma}$ is closed, whereupon it follows that the
algebra $L^{\pi}$ is closed. If $L^{\pi}$ contains a non-zero finite rank
operator $\sum_{\gamma}a_{\gamma}\otimes\pi({\gamma})$ then since each
$a_{\gamma}$ must be an operator of finite rank, the hypothesis implies that
each $a_{\gamma}=0$, and hence that $a=0$. It follows that $L^{\pi}$ does not
contain any nonzero finite rank operators, and hence that it is Engel by
Theorem \ref{engel}. The result now follows by applying Lemma \ref{maptri}.
\end{proof}

For a $\Gamma$-subgraded Lie algebra $L$, define
\[
\mathcal{C}_{\Gamma}(L)=\{[a,b]:a,b\in\cup_{\gamma\in\Gamma}L_{\gamma}\}.
\]
Then $\mathcal{C}_{\Gamma}(L)$ is the set of all \textit{homogeneous}
\textit{commutators} in $L$. It is obvious that $\mathcal{C}_{\Gamma}(L)$ is a
Lie submultiplicative set in $L$.

Now we may apply Theorem \ref{cons} to obtain the following result.

\begin{theorem}
\label{multiset}If $L$ is a $\Gamma$-subgraded Lie algebra of compact
operators on $\mathcal{X}$, and if $\mathcal{C}_{\Gamma}(L)$ consists of Engel
elements of $L$, then $L$ is $E$-solvable and is therefore triangularizable.
\end{theorem}

\begin{proof}
It is sufficient to show that $L$ is reducible. One may assume that every
component of $L$ is a closed subspace in $\mathcal{B}(\mathcal{X})$. Let
\[
I=\mathcal{C}_{\Gamma}(L)\cap\mathcal{F}(\mathcal{X}).
\]
Then $I$ is a Lie multiplicative set in $L$.

Suppose, to the contrary, that $L$ is irreducible. Then $L$ is not Engel, and
then there is a component $L_{\gamma}$ containing a non-zero finite rank
operator, say $a$ by Theorem \ref{cons}. If $a$ is scalar then $\mathcal{X}$
is finite dimensional and $\operatorname*{span}I$ consists of Engel elements
of $L$ by Lemma \ref{multi}. But in this case $I=\mathcal{C}_{\Gamma}(L)$,
whence $\operatorname*{span}I=[L,L]$, and $L$ is triangularizable by Theorem
\ref{redcri}, a contradiction.

So one can assume that $a$ is not scalar. If $[a,L_{\beta}]=0$ for every
$\beta\in\Gamma$ then $[a,L]=0$, whence $L_{0}$ is reducible by Lomonosov's
Theorem \cite{Lom}, which gives a contradiction.

Therefore $[a,L_{\beta}]\neq0$ for some $\beta\in\Gamma$, whence
\[
I\neq0
\]
and, moreover, $I$ is non-scalar. Since $\mathcal{C}_{\Gamma}(L)$ consists of
Engel elements of $L$ by hypothesis, so is the set $I$. By Lemma \ref{multi},
$\operatorname*{span}I$ consists of Engel elements of $L$. In particular,
$\operatorname*{span}I$ is Engel. But $\operatorname*{span}I$ is an ideal of
$L$, whence $L$ is reducible by Theorem \ref{solid}, a contradiction.
\end{proof}

One consequence of Theorem \ref{multiset} is the following homogeneous version
of the local triangularizability result from \cite[Corollary 5.17]{ShT2000}.

\begin{corollary}
\label{loctriang} Let $L=\sum L_{\gamma}$ be a $\Gamma$-subgraded Lie algebra
of compact operators. If for any homogeneous $a,b$ the set $\left\{
a,b\right\}  $ is triangularizable, then $L$ is triangularizable, and hence is
$E$-solvable.
\end{corollary}

\begin{proof}
It follows from homogeneous local triangularizability of $L$ that
$\mathcal{C}_{\Gamma}(L)$ is Volterra, whence $L$ is triangularizable by
Theorem \ref{multiset}.
\end{proof}

A stronger result is contained in Corollary \ref{tria}.

\section{\label{_6finite}The Case that the Underlying Space is Finite
Dimensional}

Theorem \ref{cons} shows that one should analyze graded Lie algebras with
non-zero finite rank operators. In this section, we consider the case when $L$
is a Lie algebra of operators on a finite dimensional vector space.

We begin with a graded version of Cartan's criterion.

\begin{lemma}
\label{cart}Let $L$ be a $\Gamma$-graded Lie subalgebra of operators on a
finite dimensional space. If $L_{0}$ is scalar, and some component $L_{\alpha
}$ contains an Engel element $a$ of $L$, then $a$ belongs to the radical $R$
of $L$. If in particular, $a$ is non-scalar, then $L$ is reducible.
\end{lemma}

\begin{proof}
Using the non-graded version of Cartan's Criterion, it suffices to show that
$\left\langle x,a\right\rangle =0$ for every $x$ in $L$. For arbitrary $x$ in
$L_{\gamma}$, the operator $(\operatorname*{ad}x)(\operatorname*{ad}a)$ maps
each $L_{\delta}$ to $L_{\alpha+\gamma+\delta}$. If $\alpha+\gamma\not =0$,
then $(\operatorname*{ad}x)(\operatorname*{ad}a)$ doesn't fix any component of
$L$. Thus, an easy calculation shows that in some basis of $L$, the operator
$(\operatorname*{ad}x)(\operatorname*{ad}a)$ is represented by a matrix with a
zero main diagonal, whence
\[
\left\langle x,a\right\rangle =0.
\]
On the other hand, if $\alpha+\gamma=0$, then $[x,a]$ is in $L_{0}$. Since
$L_{0}$ is scalar, this implies that $x$ and $a$ commute, whence
$\operatorname*{ad}x$ and $\operatorname*{ad}a$ commute. Since $a$ is an Engel
element, $\operatorname*{ad}a$ is nilpotent, and it follows that
$(\operatorname*{ad}x)(\operatorname*{ad}a)$ is nilpotent, whence
$\left\langle x,a\right\rangle =0$. Hence
\[
\left\langle x,a\right\rangle =0
\]
for all $x$ in $L$. Let $I=\{u\in L:$ $\left\langle u,L\right\rangle =0\}$.
Since $[v,L]\subset L$ for every $v\in L$, it follows from
\[
\left\langle \lbrack u,v],L\right\rangle =\left\langle u,[v,L]\right\rangle
=0
\]
for every $u\in I$ that in fact the subspace $I$ is an ideal of $L$. By Lemma
\ref{cartan}, $I$ is solvable, and therefore%
\[
a\in I\subset R.
\]
If $a$ is not scalar, then $L$ is reducible by Theorem \ref{solid}.
\end{proof}

One can readily deduce the following lemma from algebraic results of
\cite{Kre, KrK}, but we prefer to give a self-contained and comparatively
simple proof for the convenience of readers.

\begin{lemma}
\label{prime}Let $L$ be a $\mathbb{Z}_{n}$-graded Lie subalgebra of operators
on a finite dimensional space $\mathcal{X}$. If $L_{0}$ is scalar, then $L$ is
solvable, and therefore is triangularizable.
\end{lemma}

\begin{proof}
One may assume that $\dim\mathcal{X}>1$. Then, proceeding inductively, it
suffices to prove that $L$ is reducible.

Suppose, to the contrary, that $L$ is irreducible. It is clear that $L$ is not scalar.

\begin{claim}
\label{prime1} If $[L_{m},[L_{m},L_{k}]]$ is scalar then $[L_{m},L_{k}]=0$.
\end{claim}

Suppose that $[a,b]\neq0$ for some $a\in L_{m}$ and $b\in L_{k}$. By
condition, $[a,[a,b]]$ is scalar. Since $\operatorname*{tr}([a,[a,b]])=0$, we
obtain that
\[
\lbrack a,[a,b]]=0.
\]
Then by the Kleinecke-Shirokov Theorem \cite{Kle, Shi}, $[a,b]$ is nilpotent,
meaning $L_{m+k}$ contains a non-scalar Engel element of $L$, whence $L$ is
reducible by Lemma \ref{cart}, which gives a contradiction.

Therefore we have that
\[
\lbrack a,b]=0,
\]
and this proves Claim \ref{prime1}.

\begin{claim}
\label{prime0} If $L_{m}$ is not scalar then there is a positive integer $k<n$
such that $[L_{m},L_{k}]$ and $[L_{m},[L_{m},L_{k}]]$ are not scalar.
\end{claim}

If $L_{m}$ is in the center of $L$ then $L$ is reducible by \cite{Lom}, a
contradiction. Thus there is an integer $k>0$ such that $k<n$, and%
\[
\lbrack L_{m},L_{k}]\neq0.
\]
It follows from Claim \ref{prime1} that $[L_{m},[L_{m},L_{k}]]$ and
$[L_{m},L_{k}]$ are non-scalar.

\begin{claim}
\label{prime2} If $L_{m}$ is not scalar and an integer $k$ is as in Claim
$\ref{prime0}$ then $L_{k+pm}$ is not scalar for every integer $p>0$.
\end{claim}

Inductively applying Claims \ref{prime1} and \ref{prime0}, we obtain that%
\[
\overset{\quad\quad\quad\quad\quad p\text{ times}}{L_{k+pm}\supset
\overbrace{[L_{m},\ldots,[L_{m},[L_{m},[L_{m},L_{k}]]]\cdots]}}%
\]
is non-scalar for every $p$.

Let $\gcd(m,n)$ be the greatest common divisor of integers $m$ and $n$.

\begin{claim}
\label{prime3} If $\gcd(m,n)=1$ for some integer $m>0$, then $L_{m}$ is scalar.
\end{claim}

Suppose, to the contrary, that $L_{m}$ is non-scalar. Then it follows from
Claim \ref{prime2} that there is an integer $k$ such that $L_{k+pm}$ is
non-scalar for every integer $p>0$. On the other hand, since $\gcd(m,n)=1$,
one can find a positive integer $p$ such that $n$ divides $k+pm$. Since
$L_{k+pm}=L_{0}$, $L_{k+pm}$ is scalar, a contradiction. This proves Claim
\ref{prime3}.

\begin{claim}
\label{prime4} $L_{n-1}$, $L_{n}$ and $L_{n+1}$ are scalar.
\end{claim}

Indeed, it follows from Claim \ref{prime3} that $L_{n-1}$ and $L_{n+1}$
($=L_{1}$) consist of scalar operators. Recall that $L_{n}$ ($=L_{0}$) is also
scalar, by hypothesis.

\begin{claim}
\label{prime5} $L_{m}$ is scalar for every positive integer $m<n$.
\end{claim}

If not, let $m$ be the least positive integer such that $L_{0},L_{1}%
,...,L_{m-1}$ consist of scalar operators, but $L_{m}$ is non-scalar. Then
there is an integer $k$ satisfying Claim \ref{prime0}. On the other hand,
there is an integer $p>0$ such that
\[
k+pm\in\{n-1,n,n+1,...,n+m-1\}.
\]
Taking into account Claim \ref{prime4}, we conclude that for any integer $q$
in this set, $L_{q}$ must be scalar. But $L_{k+pm}$ is non-scalar by Claim
\ref{prime2}, a contradiction. This proves Claim \ref{prime5}.

Now it readily follows by Claim \ref{prime5} that $L$ consists of scalar
operators, a contradiction. We have therefore proved that $L$ is
triangularizable, and hence is solvable.
\end{proof}

The following example shows that the restriction to cyclic groups in the above
lemma cannot be removed.

\begin{example}
\label{pauli}Define
\[
a=%
\begin{pmatrix}
0 & 1\\
-1 & 0
\end{pmatrix}
,\quad b=%
\begin{pmatrix}
0 & -i\\
-i & 0
\end{pmatrix}
,\quad c=%
\begin{pmatrix}
-i & 0\\
0 & i
\end{pmatrix}
.
\]
These matrices form a basis of $L:=\mathfrak{sl}(2,\mathbb{C})$, the Lie
algebra of all operators on $\mathbb{C}^{2}$ with zero trace, and satisfy the
identities
\[
\lbrack a,b]=2c,\quad\lbrack b,c]=2a,\quad\lbrack c,a]=2b.
\]
Let $\Gamma=\mathbb{Z}_{2}\times\mathbb{Z}_{2}$; then we may define a $\Gamma
$-grading on $L$ by setting
\[
L_{(0,0)}=0,\quad L_{(0,1)}=\mathbb{C}a,\quad L_{(1,0)}=\mathbb{C}b,\quad
L_{(1,1)}=\mathbb{C}c.
\]
It is clear that $L$ is irreducible.

\emph{Note that the spectral radii of} $a,b,c$ \emph{are equal to} $1$.
\emph{It is easy to calculate that}
\[
\rho(\lambda x+\mu y)=\left\vert \lambda^{2}+\mu^{2}\right\vert ^{1/2}%
\leq\left\vert \lambda\right\vert +\left\vert \mu\right\vert =\rho(\lambda
x)+\rho(\mu y),
\]
\emph{and}%
\[
\rho(\lambda x\mu y)=\left\vert \lambda\right\vert \left\vert \mu\right\vert
\rho(xy)=\left\vert \lambda\right\vert \left\vert \mu\right\vert =\rho(\lambda
x)\rho(\mu y),
\]
\emph{for every} $\lambda,\mu\in\mathbb{C}$, \emph{and} $x,y\in\{a,b,c\}$
\emph{with} $x\neq y$. \emph{This shows that the following conditions (even
taken together)}

\begin{itemize}
\item $\rho(u+v)\leq\rho(u)+\rho(v)$ for every $u,v\in\cup_{\gamma\in\Gamma
}L_{\gamma}$.

\item $\rho(uv)\leq\rho(u)\rho(v)$ for every $u,v\in\cup_{\gamma\in\Gamma
}L_{\gamma}$.
\end{itemize}

\noindent\emph{are not sufficient to imply the reducibility of a} $\Gamma
$-\emph{graded Lie algebra of compact operators with }$L_{0}=0$.
\end{example}

\begin{lemma}
\label{finsubgraded} Let $L$ be a $\Gamma$-subgraded Lie algebra of operators
on a finite dimensional space. Then $L$ is solvable if one of the following
conditions holds.

\begin{itemize}
\item[(i)] $L_{\alpha}$ consists of Engel elements of $L$ for every $\alpha
\in\Gamma$.

\item[(ii)] $\Gamma=\mathbb{Z}_{n}$ and $L_{0}$ consists of Engel elements of
$L$.
\end{itemize}
\end{lemma}

\begin{proof}
First, assume that $L$ is graded. By Theorem \ref{redcri}, we must show that
$L$ is triangularizable. For this, it suffices to prove that $L$ is reducible.

\begin{case}
$L_{0}$ is not scalar.
\end{case}

For $\alpha\not =0$, as in Lemma \ref{cart}, we obtain that
\[
\left\langle L_{0},L_{\alpha}\right\rangle =0.
\]
Since $\operatorname*{ad}L_{0}$ consists of nilpotent elements, we also have
\[
\left\langle L_{0},L_{0}\right\rangle =0.
\]
Since the Killing form is linear, we therefore obtain that
\[
\left\langle L_{0},L\right\rangle =0.
\]
Let $I=\{a\in L:\left\langle a,L\right\rangle =0\}$. Then $I$ is a solvable
ideal by Lemma \ref{cartan}. Since
\[
L_{0}\subset I,
\]
$I$ is non-scalar, whence $L$ is reducible by Theorem \ref{solid}. This proves
both (i) and (ii).

\begin{case}
$L_{0}$ is scalar.
\end{case}

In this case, (i) follows from Lemma \ref{cart}, and (ii) follows from Lemma
\ref{prime}.

\begin{case}
$L$ is not necessarily graded (but of course is $\Gamma$-subgraded).
\end{case}

Let $M=L^{\pi}$ (see the beginning of Section \ref{_5tensor}). Then $M$ is a
$\Gamma$-graded Lie algebra of compact operators on a finite dimensional
space. It is easy to see that if some $L_{\gamma}$ consists of Engel elements
of $L$ then $M_{\gamma}$ consists of Engel elements of $M$. Indeed, one can
check that
\[
\operatorname*{ad}\nolimits_{M}(a\otimes\pi(\gamma))=(\operatorname*{ad}%
\nolimits_{L}a)\otimes\pi^{+}(\gamma)
\]
where $\pi^{+}$ is the representation of $\Gamma$ on the image of $\pi$
defined by
\[
\pi^{+}(\gamma)\pi(\alpha)=\pi(\gamma+\alpha)
\]
for every $\alpha\in\Gamma$. Now if $\operatorname*{ad}_{L}a$ is nilpotent,
then so is $(\operatorname*{ad}\nolimits_{L}a)\otimes\pi^{+}(\gamma)$, and we
conclude that $\operatorname*{ad}\nolimits_{M}(a\otimes\pi(\gamma))$ must also
be nilpotent. It now follows from above that $M$ is solvable; the solvability
of $L$ follows from Lemma \ref{maptri}.
\end{proof}

In the section on $\mathbb{Z}_{2}$-graded algebras we will make use of the
following result which is reminiscent of Cartan's Criterion.

\begin{lemma}
\label{L0-triang} Let $L$ be a $\Gamma$-graded Lie algebra of operators on a
finite dimensional space. If $L_{0}$ is solvable and non-commutative, then $L$
has a non-scalar solvable ideal, and is therefore reducible.
\end{lemma}

\begin{proof}
By Lemma \ref{cartan}, we have that
\[
\left\langle \lbrack L_{0},L_{0}],L_{0}\right\rangle =0.
\]
Since, as in Lemma \ref{cart}, $\left\langle L_{0},L_{\alpha}\right\rangle =0$
for all $\alpha\neq0$, and $[L_{0},L_{0}]\subset L_{0}$, we obtain that
\[
\left\langle \lbrack L_{0},L_{0}],L_{\alpha}\right\rangle =0
\]
for all $\alpha\neq0$. This implies that
\[
\left\langle \lbrack L_{0},L_{0}],L\right\rangle =0.
\]
Let $I=\{a\in L:\left\langle a,L\right\rangle =0\}$. As we have seen, $I$ is a
solvable ideal of $L$ by Lemma \ref{cartan}. Note that it is non-scalar
because
\[
\lbrack L_{0},L_{0}]\subset I,
\]
and $[L_{0},L_{0}]$ contains a non-zero nilpotent operator. Therefore, $L$ is
reducible by Theorem \ref{solid}.
\end{proof}

\section{ \label{_7triangularization}Triangularization Theorems}

Taking Theorem \ref{cons} into account, we have to consider those subgraded
Lie algebras of compact operators on infinite dimensional spaces whose
components contain finite rank operators. Recall (Remark \ref{adnil}) that a
finite rank operator $a$ in an operator Lie algebra $L$ is an Engel element of
$L$ if and only if the operator $\operatorname*{ad}a$ is nilpotent on $L$.

\begin{lemma}
\label{finrank} Suppose that a $\Gamma$-subgraded Lie algebra $L$ consists of
finite rank operators and that all components $L_{\gamma}$ consist of Engel
elements of $L$ (if $\Gamma=\mathbb{Z}_{n}$, then we impose the last condition
only on operators in $L_{0}$). Then $L$ is $E$-solvable and is therefore triangularizable.
\end{lemma}

\begin{proof}
Given a finite dimensional subspace $\mathcal{Y}$ of $\mathcal{X}$, let
\[
L_{\gamma}^{\mathcal{Y}}=\{a\in L_{\gamma}:a\mathcal{X}\subset\mathcal{Y}%
\}\text{ and }L^{\mathcal{Y}}=\sum_{\gamma\in\Gamma}L_{\gamma}^{\mathcal{Y}}.
\]
Then $L^{\mathcal{Y}}$ is a $\Gamma$-subgraded Lie algebra of operators on
$\mathcal{X}$. It is evident that $\mathcal{Y}$ is invariant for
$L^{\mathcal{Y}}$, and the restriction $L^{\mathcal{Y}}|\mathcal{Y}$ of
$L^{\mathcal{Y}}$ to $\mathcal{Y}$ is a $\Gamma$-subgraded Lie algebra of
operators on $\mathcal{Y}$, namely
\[
L^{\mathcal{Y}}|\mathcal{Y}=\sum_{\gamma\in\Gamma}L_{\gamma}^{\mathcal{Y}%
}|\mathcal{Y}.
\]
Since $\mathcal{Y}$ is finite dimensional, $L^{\mathcal{Y}}|\mathcal{Y}$ is
also finite dimensional.

Suppose that some $L_{\gamma}$ consists of Engel elements of $L$, and let $a$
be an arbitrary element of $L_{\gamma}^{\mathcal{Y}}$. Since $a\in L_{\gamma}%
$, we obtain for every $b\in L^{\mathcal{Y}}$ that
\[
\left\Vert \left(  \operatorname*{ad}a|\mathcal{Y}\right)  ^{n}\left(
b|\mathcal{Y}\right)  \right\Vert ^{1/n}=\left\Vert (\left(
\operatorname*{ad}a\right)  ^{n}b)|\mathcal{Y}\right\Vert ^{1/n}\leq\left\Vert
\left(  \operatorname*{ad}a\right)  ^{n}b\right\Vert ^{1/n}\rightarrow0
\]
as $n\rightarrow\infty$, whence $\operatorname*{ad}a|\mathcal{Y}$ is nilpotent
on $L^{\mathcal{Y}}|\mathcal{Y}$. Therefore, the component $L_{\gamma
}^{\mathcal{Y}}|\mathcal{Y}$ of $L^{\mathcal{Y}}|\mathcal{Y}$ consists of
Engel elements of $L^{\mathcal{Y}}|\mathcal{Y}$. Thus $L^{\mathcal{Y}%
}|\mathcal{Y}$ satisfies conditions of Lemma \ref{finsubgraded}. Then, by
Lemma \ref{finsubgraded}, $L^{\mathcal{Y}}|\mathcal{Y}$ is solvable and is
triangularizable, whence any sum of commutators in $L^{\mathcal{Y}%
}|\mathcal{Y}$ is nilpotent.

Let $a_{i},b_{i}$ ($1\leq i\leq n$) be arbitrary operators in $L$. Then
$a_{i}=\sum_{\gamma}a_{i,\gamma}$ and $b_{i}=\sum_{\gamma}b_{i,\gamma}$, where
the sums are finite, in accordance of grading on $L$. This implies the
existence of a finite dimensional subspace $\mathcal{Y}$ that contains all
$a_{i,\gamma}\mathcal{X}$ and $b_{i,\gamma}\mathcal{X}$. It follows that
$a_{i},b_{i}\in L^{\mathcal{Y}}$ for every $i$. Let
\[
c=\sum_{i}[a_{i},b_{i}].
\]
If $c$ is not nilpotent then there is a non-zero vector $x\in\mathcal{X}$ such
that
\[
cx=\lambda x
\]
for some non-zero $\lambda\in\mathbb{C}$; but $x\in\mathcal{Y}$ and
$c|\mathcal{Y}$ is not nilpotent, which gives a contradiction.

This proves that $[L,L]$ consists of nilpotent operators, and hence that $L$
is $E$-solvable and triangularizable by Theorem \ref{redcri}.
\end{proof}

\begin{lemma}
\label{onefin} Let $L=\sum L_{\gamma}$ be a $\Gamma$-subgraded Lie algebra of
compact operators. Suppose that $\cup_{\gamma\in\Gamma}L_{\gamma}$ contains a
non-scalar finite rank operator and each $L_{\gamma}$ consists of Engel
elements of $L$ (if $\Gamma=\mathbb{Z}_{n}$, then we again impose this
condition only on $L_{0}$). Then $L$ contains a non-scalar $E$-solvable ideal.
\end{lemma}

\begin{proof}
Let $I_{\gamma}=L_{\gamma}\cap\mathcal{F}(\mathcal{X})$ and $I=\sum_{\gamma
\in\Gamma}I_{\gamma}$. Then $I$ is a non-scalar ideal of $L$, so it suffices
to show that $I$ is triangularizable. But this follows from Lemma
\ref{finrank}, because $I$ is a $\Gamma$-subgraded Lie algebra of finite rank
operators and each $I_{\gamma}$ consists of Engel elements of $I$ (if
$\Gamma=\mathbb{Z}_{n}$, then this need be true only of $I_{0}$).
\end{proof}

Now we obtain our main result on triangularization of $\Gamma$-subgraded Lie algebras.

\begin{theorem}
\label{main1} Let $L$ be a $\Gamma$-subgraded (respectively, $\mathbb{Z}_{n}%
$-subgraded) Lie algebra of compact operators. If each $L_{\gamma}$
(respectively, only $L_{0}$) consists of Engel elements of $L$, then $L$ is
$E$-solvable, and is therefore triangularizable.
\end{theorem}

\begin{proof}
Follows from Lemma \ref{onefin} and Theorem \ref{cons}.
\end{proof}

\begin{corollary}
\label{volt} If the components of a $\Gamma$-subgraded Lie algebra $L$ of
compact operators consist of Engel elements of $L$ (respectively, of Volterra
operators) then $L$ is Engel (respectively, Volterra).
\end{corollary}

\begin{proof}
Indeed, $L$ is triangularizable by Theorem \ref{main1}. Since the spectral
radius is subadditive on $\operatorname*{ad}L$ (respectively, on $L$) by
Theorem \ref{redcri}, and since every operator in $L$ is a sum of elements of
components of $L$, every operator in $L$ is an Engel element of $L$
(respectively, is a Volterra operator).
\end{proof}

The following example shows that one cannot weaken the assumptions of Theorem
\ref{main1} to the case that the components are triangularizable.

\begin{example}
\label{e1}Let $\mathcal{X}$ be two-dimensional and $L=\mathfrak{sl}%
(2,\mathbb{C})$. Let $e,f,g$ be operators in $L$ such that
\[
\lbrack e,f]=g,\quad\lbrack g,e]=e,\quad\lbrack g,f]=-f.
\]
Then $L$ has a $\mathbb{Z}_{3}$-grading, defined by setting
\[
L_{0}=\mathbb{C}g,\quad L_{1}=\mathbb{C}e,\quad L_{2}=\mathbb{C}f.
\]
The components $L_{1}$ and $L_{2}$ are Volterra, $L_{0}$ is triangularizable
(more precisely, is commutative), but $L$ is clearly irreducible.

\emph{Note that }$L_{0}$ \emph{doesn't consist of Engel elements of} $L$,
\emph{and that} $[L_{1},L_{2}]$ \emph{doesn't consist of Volterra operators.
The commutativity of} $L_{0}$ \emph{is a necessary condition for such examples
in virtue of Lemma} $\ref{L0-triang}$. \emph{It follows from Remark}
$\ref{z3}$ \emph{below that the condition} $[L_{1},L_{2}]$ \emph{is Volterra
is essential for reducibility of} $\mathbb{Z}_{3}$-\emph{subgraded Lie
algebras of compact operators with triangularizable} $L_{0}$.
\end{example}

For a $\Gamma$-subgraded Lie algebra $L$, let
\[
L_{0}^{\prime}=\sum_{\gamma\in\Gamma}[L_{\gamma},L_{-\gamma}]\text{\quad
and\quad}L_{0}^{\prime\prime}=\sum_{\gamma\in\Gamma\backslash\{0\}}[L_{\gamma
},L_{-\gamma}].
\]
It is clear that
\[
L_{0}^{\prime\prime}\subset L_{0}^{\prime}\subset L_{0}.
\]
We also set
\[
L^{\prime}=L_{0}^{\prime}+\sum_{\gamma\in\Gamma\backslash\{0\}}L_{\gamma
}\text{\quad and\quad}L^{\prime\prime}=L_{0}^{\prime\prime}+\sum_{\gamma
\in\Gamma\backslash\{0\}}L_{\gamma}.
\]
It is easy to see that $L^{\prime}$ and $L^{\prime\prime}$ are $\Gamma
$-subgraded Lie algebras whose zero components are $L_{0}^{\prime}$ and
$L_{0}^{\prime\prime}$, respectively, and that
\[
L_{\gamma}^{\prime}=L_{\gamma}^{\prime\prime}=L_{\gamma}%
\]
for every $\gamma\in\Gamma\backslash\{0\}$.

The following lemma motivates the definition of $L^{\prime}$ and
$L^{\prime\prime}$.

\begin{lemma}
\label{ideal} For a $\Gamma$-subgraded Lie algebra $L$ of compact operators,
the Lie algebras $L^{\prime}$ and $L^{\prime\prime}$ are ideals of $L$.
Moreover, $L$ is $E$-solvable if and only if is $L^{\prime}$.
\end{lemma}

\begin{proof}
It is easy to check that $[L^{\prime},L]\subset L^{\prime}$ and $[L^{\prime
\prime},L]\subset L^{\prime\prime}$.

If $L$ is $E$-solvable, then $L$ is triangularizable, whence $L^{\prime}$ is
triangularizable, and is therefore $E$-solvable.

If $L^{\prime}$ is $E$-solvable, then taking into account that
\[
\lbrack L,L]\subset L^{\prime},
\]
we have that $[L,L]$ is triangularizable, whence $L$ is triangularizable, and
is therefore $E$-solvable.
\end{proof}

Define the relation $\Gamma^{\sharp}$ in $\Gamma\times\Gamma$ by setting%
\[
\Gamma^{\sharp}=\{(\alpha,\beta)\in\Gamma\times\Gamma:\text{ no cyclic
subgroup of }\Gamma\text{ contains }\alpha\text{ and }\beta\}.
\]
It is clear that $\Gamma^{\sharp}$ is an empty relation if $\Gamma$ is cyclic.

Now we are in a position to obtain a stronger version of Theorem
\ref{multiset}.

\begin{theorem}
\label{lieset}Let $L$ be a $\Gamma$-subgraded Lie algebra of compact
operators. If every commutator $[a,b]$ is an Engel element of $L$ for $a\in
L_{\gamma}$ and $b\in L_{\delta}$ such that $\gamma+\delta=0$ or
$(\gamma,\delta)\in\Gamma^{\sharp}$, then $L$ is $E$-solvable and is triangularizable.
\end{theorem}

\begin{proof}
Since $\mathcal{C}(L_{0})$ consists of Engel elements of $L$, it follows from
Theorem \ref{crit} that $L_{0}$ is triangularizable.

\begin{claim}
\label{lieset2} $L_{0}^{\prime}$ consists of Engel elements of $L^{\prime}$.
\end{claim}

Indeed, since
\[
L_{0}^{\prime}\subset L_{0},
\]
we have that $L_{0}^{\prime}$ is triangularizable. By hypothesis, every
commutator $[a,b]$, with $a\in L_{\gamma}$ and $b\in L_{-\gamma}$, is an Engel
element of $L$ for every $\gamma\in\Gamma$. Since a finite sum of Engel
elements of $L$ in the triangularizable Lie algebra $L_{0}$ is again an Engel
element of $L$ by Theorem \ref{tensor}, it follows that $L_{0}^{\prime}$
consists of Engel elements of $L$. Since
\[
L^{\prime}\subset L,
\]
$L_{0}^{\prime}$ consists of Engel elements of $L^{\prime}$.

\begin{claim}
\label{lieset3} The set $\mathcal{C}_{\Gamma}(L^{\prime})$ of homogeneous
commutators in $L^{\prime}$ consists of Engel elements of $L^{\prime}$.
\end{claim}

Let $a\in L_{\gamma}^{\prime}$ and $b\in L_{\delta}^{\prime}$. Assume first
that there is a cyclic subgroup $\Lambda\subset\Gamma$ containing $\gamma$ and
$\delta$. Then
\[
a,b\in N:=\sum_{\beta\in\Lambda}L_{\beta}^{\prime},
\]
where $N$ is a $\Lambda$-subgraded Lie algebra of compact operators. Since the
zero component of $N$ consists of Engel elements of $L^{\prime}$ by Claim
\ref{lieset2}, this component consists of Engel elements of $N$. Then $N$ is
triangularizable by Theorem \ref{main1}, whence $[a,b]$ is a Volterra
operator. In particular, $[a,b]$ is an Engel element of $L^{\prime}$.

If $(\gamma,\delta)\in\Gamma^{\sharp}$ then $[a,b]$ is an Engel element of
$L$, and hence of $L^{\prime}$, by hypothesis.

\begin{claim}
\label{lieset4} $L$ is $E$-solvable and triangularizable.
\end{claim}

Indeed, by Claim \ref{lieset3} and Theorem \ref{multiset}, $L^{\prime}$ is
triangularizable. By Theorem \ref{redcri}, $L^{\prime}$ is $E$-solvable. Now
the claim follows from Lemma \ref{ideal}.
\end{proof}

Theorem \ref{lieset} allows us to obtain the following result which extends
Corollary \ref{loctriang}.

\begin{corollary}
\label{tria}Let $L$ be a $\Gamma$-subgraded Lie algebra of compact operators.
If every pair $\{a,b\}$ is triangularizable for $a\in L_{\gamma}$ and $b\in
L_{\delta}$ such that $\gamma+\delta=0$ or $(\gamma,\delta)\in\Gamma^{\sharp}%
$, then $L$ is triangularizable.
\end{corollary}

\begin{proof}
Indeed, every commutator pointed out in Theorem $\ref{lieset}$ is Volterra in
$L$, whence $L$ is $E$-solvable and, by Theorem $\ref{redcri}$, is triangularizable.
\end{proof}

\section{\label{_8consequences}Consequences}

\begin{corollary}
\label{main2} Let $L$ be a $\Gamma$-subgraded Lie algebra of compact
operators. Then $L$ is $E$-solvable if one of the following conditions holds.

\begin{itemize}
\item[(i)] $L_{0}$ consists of Engel elements of $L$, and $[L_{\gamma
},L_{\delta}]$ consists of Engel elements of $L$ for every non-zero
$\gamma,\delta\in\Gamma$ with $(\gamma,\delta)\in\Gamma^{\sharp}$.

\item[(ii)] $L_{0}$ is $E$-solvable, and $[L_{\gamma},L_{\delta}]$ consists of
Engel elements of $L$ for every non-zero $\gamma,\delta\in\Gamma$ such that
$\gamma+\delta=0$ or $(\gamma,\delta)\in\Gamma^{\sharp}$.
\end{itemize}
\end{corollary}

\begin{proof}
First, we check that the conditions of Theorem \ref{lieset} hold in each case.

For (ii), it is sufficient to show that $[L_{0},L_{0}]$ is Volterra. Since
$L_{0}$ is $E$-solvable, this is obvious by Theorem \ref{redcri}.

For (i), it is sufficient to show that $[L_{\gamma},L_{-\gamma}]$ is Volterra
for every $\gamma\in\Gamma$. For $\gamma=0$, we have that $[L_{0},L_{0}]$ is
Volterra, because $L_{0}$ is triangularizable. Now take a non-zero $\gamma
\in\Gamma$, and let $\Lambda$ be a subgroup of $\Gamma$ generated by $\gamma$.
Then $\Lambda$ is cyclic. Let
\[
N=\sum_{\alpha\in\Lambda}L_{\alpha}.
\]
Then $N$ is a $\Lambda$-subgraded Lie algebra of compact operators. Since
$\Lambda$ is cyclic, and since the component $L_{0}$ consists of Engel
elements of $L$, and hence of $N$, this Lie algebra $N$ is triangularizable by
Theorem \ref{main1}. Then $\left[  N,N\right]  $ is Volterra by Theorem
\ref{redcri}, and since
\[
\lbrack L_{\gamma},L_{-\gamma}]\subset\left[  N,N\right]  ,
\]
we have that $[L_{\gamma},L_{-\gamma}]$ is Volterra.

Therefore, $L$ is $E$-solvable in both cases by Theorem \ref{lieset}.
\end{proof}

In particular, Corollary \ref{main2} implies the following result.

\begin{corollary}
\label{main3} Let $L$ be a $\mathbb{Z}_{n}$-subgraded Lie algebra of compact
operators. If $L_{0}$ is $E$-solvable, and $\left[  L_{k},L_{-k}\right]  $
consists of Engel elements of $L$ for every integer $k>0$ such that $k\leq
n/2$, then $L$ is $E$-solvable and is therefore triangularizable.
\end{corollary}

\begin{proof}
Indeed, $L$ satisfies conditions of Corollary \ref{main2}(ii).
\end{proof}

\begin{remark}
\label{z3}$\emph{The}$ $\emph{special}$\emph{ }$\emph{case}$\emph{ }%
$\emph{of}$\emph{ }$\emph{Corollary}$ $\ref{main3}$ $\emph{for}$
$\mathbb{Z}_{3}$\textrm{-}$\emph{subgraded}$ $\emph{Lie}$ $\emph{algebras}$
$\emph{gives}$ $\emph{us}$ $\emph{the}$ $\emph{following}$ $\emph{result:}%
$\medskip

Let $L=L_{0}+L_{1}+L_{2}$ be a $\mathbb{Z}_{3}$-subgraded Lie algebra of
compact operators. If $L_{0}$ is $E$-solvable and $\left[  L_{1},L_{2}\right]
$ consists of Volterra operators then $L$ is $E$-solvable.\medskip

$\emph{The}$ $\emph{condition}$ $\emph{that}$ $\left[  L_{1},L_{2}\right]  $
$\emph{is}$ $\emph{Volterra}$ $\emph{is}$ $\emph{essential}$ $\emph{in}$
$\emph{virtue}$ $\emph{of}$ $\emph{Example}$ $\ref{e1}$.
\end{remark}

Both of the assertions of Theorem \ref{main1} are contained in the following result.

\begin{corollary}
\label{subgroup}Let $L$ be a $\Gamma$-subgraded Lie algebra of compact
operators. Suppose that $\Gamma$ has a subgroup $\Lambda$, such that
$\Gamma/\Lambda$ is cyclic. If $L_{\delta}$ consists of Engel elements of $L$
for every $\delta\in\Lambda$, then $L$ is $E$-solvable.
\end{corollary}

\begin{proof}
For $\alpha\in\Gamma/\Lambda$, let $M_{\alpha}=\sum_{\gamma\in\alpha}%
L_{\gamma}$. Then $M:=\sum_{\alpha\in\Gamma/\Lambda}M_{\alpha}$ is a
$(\Gamma/\Lambda)$-subgraded Lie algebra of compact operators and $M=L$. Note
that
\[
M_{0}=\sum_{\delta\in\Lambda}L_{\delta}.
\]
By Corollary \ref{volt}, $M_{0}$ is an Engel Lie algebra. Since $M_{0}$ is
triangularizable by Theorem \ref{redcri}, $M_{0}$ consists of Engel elements
of $M$ by Theorem \ref{tensor}. So $M$ is $E$-solvable by Theorem \ref{main1}.
\end{proof}

Theorem $\ref{main1}$ allows us, with certain hypotheses in place, to extend
the triangularizability of the zero component of a $\mathbb{Z}_{n}$-subgraded
Lie algebra of compact operators to the entire algebra. It is natural to ask
for the cases in which \textbf{ }one can repeat this procedure of inflation.
More specifically, \textit{let} $M$ \textit{be a} $\mathbb{Z}_{m}%
$\textit{-graded Lie algebra, and let} $N:=M_{0}$ \textit{be a} $\mathbb{Z}%
_{n}$\textit{-graded Lie algebra such that }$N_{0}$\textit{ is Volterra. Is it
true that }$M$\textit{ is }$E$\textit{-solvable?}

The following example shows that the answer to this question is negative in general.

\begin{example}
Let $L$ be as in Example $\ref{pauli}$, and define%
\[
M_{0}=L_{(0,0)}+L_{(0,1)},\quad M_{1}=L_{(1,0)}+L_{(1,1)},
\]
and%
\[
N_{0}=L_{(0,0)},\quad N_{1}=L_{(0,1)}.
\]
Then $M:=M_{0}+M_{1}$ and $N:=N_{0}+N_{1}$ are $\mathbb{Z}_{2}$-graded Lie
algebras, $N=M_{0}$, and $N_{0}=0$, but it is clear that $M=L$ is irreducible.
\end{example}

\begin{lemma}
\label{endo}Let $L$ be a normed Lie algebra (respectively, a normed algebra),
and let $\varphi$ be a bounded endomorphism of $L$. Then $\left[
\mathcal{E}_{\lambda}(\varphi),\mathcal{E}_{\mu}(\varphi)\right]
\subset\mathcal{E}_{\lambda\mu}(\varphi)$ (respectively, $\mathcal{E}%
_{\lambda}(\varphi)\mathcal{E}_{\mu}(\varphi)\subset\mathcal{E}_{\lambda\mu
}(\varphi)$).
\end{lemma}

\begin{proof}
Put $\varphi_{\lambda}=\varphi-\lambda$ for every $\lambda\in\mathbb{C}$. Let
$x,y$ be arbitrary elements of $L$.

\begin{claim}
\label{endo1} $\varphi_{\lambda\mu}(\left[  x,y\right]  )=\left[
\varphi_{\lambda}(x),\varphi_{\mu}(y)\right]  +\left[  \lambda x,\varphi_{\mu
}(y)\right]  +\left[  \varphi_{\lambda}(x),\mu y\right]  $.
\end{claim}

It is an easy calculation.

\begin{claim}
\label{endo2} $\varphi_{\lambda\mu}^{n}(\left[  x,y\right]  )$ is the sum of
$3^{n}$ summands of type $\left[  \lambda^{n-i}\varphi_{\lambda}^{i}%
(x),\mu^{n-j}\varphi_{\mu}^{j}(y)\right]  $ for integers $i$ and $j$ such that
$0\leq i\leq n$, $0\leq j\leq n$ and $i+j\geq n$.
\end{claim}

The claim follows easily by induction. Indeed, let $z_{i}=\lambda^{n-i}%
\varphi_{\lambda}^{i}(x)$ and $w_{j}=\mu^{n-j}\varphi_{\mu}^{j}(y)$ for
integers $i,j$ such that $0\leq i\leq n$, $0\leq j\leq n$ and $i+j\geq n$. It
is sufficient to show that $\varphi_{\lambda\mu}\left(  \left[  z_{i}%
,w_{j}\right]  \right)  $ may be written as the sum of three terms satisfying
the condition of Claim \ref{endo2} for $n+1$, and this follows by Claim
\ref{endo1}.

\begin{claim}
\label{endo3} $\left[  \mathcal{E}_{\lambda}(\varphi),\mathcal{E}_{\mu
}(\varphi)\right]  \subset\mathcal{E}_{\lambda\mu}(\varphi)$.
\end{claim}

Let $x\in\mathcal{E}_{\lambda}(\varphi)$ and $y\in\mathcal{E}_{\mu}(\varphi)$.
Then for every $\varepsilon>0$ there is a constant $C>0$ such that
\[
\left\Vert \varphi_{\lambda}^{m}(x)\right\Vert \leq C\varepsilon^{m}\text{ and
}\left\Vert \varphi_{\mu}^{m}(y)\right\Vert \leq C\varepsilon^{m}%
\]
for every integer $m\geq0$. Let $t=\max\left\{  \left\vert \lambda\right\vert
,\left\vert \mu\right\vert ,1\right\}  $. Assume that $\varepsilon\leq1$.
Then, by Claim \ref{endo2}, we obtain that
\begin{align*}
\left\Vert \varphi_{\lambda\mu}^{n}(\left[  x,y\right]  )\right\Vert  &
\leq3^{n}\max\left\{  \left\Vert \left[  \lambda^{n-i}\varphi_{\lambda}%
^{i}(x),\mu^{n-j}\varphi_{\mu}^{j}(y)\right]  \right\Vert :i+j\geq n\right\}
\\
&  \leq3^{n}t^{2n-i-j}2C^{2}\varepsilon^{i+j}\leq2C^{2}(3t\varepsilon)^{n}%
\end{align*}
for every integer $n\geq0$. Hence
\[
\lim\sup\left\Vert \varphi_{\lambda\mu}^{n}(\left[  x,y\right]  )\right\Vert
^{1/n}\leq3t\varepsilon.
\]
Since $\varepsilon$ is arbitrary, we have that
\[
\left[  x,y\right]  \in\mathcal{E}_{\lambda\mu}(\varphi).
\]
This proves Claim \ref{endo3}.

The case of an endomorphism of a normed algebra is similar.
\end{proof}

\begin{corollary}
\label{aut} Let $L$ be a Lie algebra of compact operators, $\varphi$ an
automorphism of $L$. Then $L$ is $E$-solvable if one of the following
conditions holds.

\begin{enumerate}
\item[(i)] $\varphi$ has a finite order, say $\varphi^{n}=1$, and the set of
fixed points of $\varphi$ consists of Engel elements of $L$.

\item[(ii)] $L$ is closed, $\varphi$ is bounded, there is an integer $n>0$
such that $\varphi^{n}-1$ is a quasinilpotent operator on $L$, and
$\mathcal{E}_{1}(\varphi)$ consists of Engel elements of $L$.
\end{enumerate}
\end{corollary}

\begin{proof}
Let $\theta=e^{2\pi i/n}$.

(i) For each $k\in\{0,1,\ldots,n-1\}$, let
\[
L_{k}=\{x\in L:\varphi(x)=\theta^{k}x\}.
\]
Then it is not difficult to check that%
\[
L=L_{0}+\ldots+L_{n-1},
\]
and that $L_{0}$ is the Lie algebra of fixed points of $\varphi$. It is also
clear that
\begin{equation}
\lbrack L_{k},L_{j}]\subset L_{k+j} \label{graded}%
\end{equation}
where the addition is modulo $n$. In other words, $L$ is $\mathbb{Z}_{n}%
$-graded. Since $L_{0}$ consists of Engel elements of $L$, $L$ is $E$-solvable
by Theorem \ref{main1}.

(ii) For each $k\in\{0,1,\ldots,n-1\}$, let
\[
L_{k}=\mathcal{E}_{\theta^{k}}(\varphi).
\]
Since $\sigma(\varphi)=\left\{  1,\theta,\ldots,\theta^{n-1}\right\}  $ and
$L_{k}$ is the image of the Riesz projection of $\varphi$ corresponding to
$\theta^{k}$, we have that%
\[
L=L_{0}+\ldots+L_{n-1},
\]
and $L_{0}=\mathcal{E}_{1}(\varphi)$. It follows from Lemma \ref{endo} that
(\ref{graded}) holds. Therefore $L$ is a $\mathbb{Z}_{n}$-graded Lie
algebra\textbf{ }and $L$ is $E$-solvable.
\end{proof}

\begin{remark}
$\emph{Recall}$ $\emph{that}$ if $L$ is $\mathbb{Z}_{n}$-graded, then there is
an automorphism $\varphi$ of $L$ such that $\varphi^{n}=1$ and $L_{0}=\{a\in
L:\varphi(a)=a\}$. \emph{Indeed, it suffices to set} $\varphi(a)=\theta^{k}a$
\emph{for} $a\in L_{k}$ \emph{and extend} $\varphi$ \emph{to} $L$ \emph{by
linearity}.

$\emph{Thus}$ $\emph{Corollary}$ $\ref{aut}(\mathrm{i})$ $\emph{is}$
$\emph{in}$ $\emph{fact}$ $\emph{a}$ $\emph{reformulation}$ $\emph{of}$
$\emph{the}$ $\emph{part}$ $\emph{of}$ $\emph{Theorem}$ $\ref{main1}%
\mathfrak{,}$ $\emph{which}$ $\emph{tells}$ $\emph{about}$ $\mathbb{Z}_{n}%
$-$\emph{graded}$ $\emph{Lie}$ $\emph{algebras}$.
\end{remark}

\begin{corollary}
\label{onetwo} Let $L$ be a $\mathbb{Z}_{n}$-subgraded Lie algebra of compact
operators, and let $I$ be $L^{\prime}$ or $L^{\prime\prime}$. If at least one
of \textbf{ }the components of $I$ is non-scalar, and $I_{0}$ consists of
Engel elements of $I$, then $L$ is reducible. In particular, if $I=L^{\prime}$
then $L$ is $E$-solvable.
\end{corollary}

\begin{proof}
Since $I$ is an ideal of $L$, it follows from Theorem \ref{main1} that $I$ is
an $E$-solvable Lie algebra, and hence $L$ is reducible by Theorem
\ref{solid}. If $I=L^{\prime}$, then it follows from Lemma \ref{ideal} that
$L$ is $E$-solvable.
\end{proof}

\begin{corollary}
Let $L$ be a $\mathbb{Z}_{n}$-subgraded Lie algebra of compact operators. If
$L_{k}$ commutes with $L_{-k}$ for every $k\neq0$, and at least one of
components of $L^{\prime\prime}$ is not scalar, then $L$ is reducible.
\end{corollary}

\begin{proof}
By hypothesis, $L_{0}^{\prime\prime}=0$, whence $L$ is reducible by Corollary
\ref{onetwo}.
\end{proof}

The following statement is a variation of Theorem \ref{cons}.

\begin{corollary}
\label{nilfinrank}Let $L$ be a $\mathbb{Z}_{n}$-subgraded Lie algebra of
compact operators, and let $I$ be one of $L^{\prime\prime}$, $L^{\prime}$ and
$L$. If at least one of the components of $I$ is non-scalar, and if there are
no non-zero nilpotent finite rank operators in the the closure of any $I_{k}$,
then $I_{0}$ consists of Engel elements of $I$, and $L$ is reducible. In
particular, if $I=L$ or $I=L^{\prime}$, then $L$ is $E$-solvable.
\end{corollary}

\begin{proof}
Let $a$ be an arbitrary element of $L_{0}$. If $(\operatorname*{ad}%
a)|\overline{I_{k}}$ has a non-zero isolated point $\lambda$ in its spectrum
for some $k$ then there is a non-zero $x\in\overline{I_{k}}$ such that
\[
x\in\mathcal{E}_{\lambda}(\operatorname*{ad}a).
\]
By \cite[Theorem 3]{Woj}, $x$ is a nilpotent finite rank operator, which would
be a contradiction.

So $\operatorname*{ad}a$ is locally quasinilpotent on every $I_{k}$ and
therefore on their sum, that is,
\[
I\subset\mathcal{E}_{0}(\operatorname*{ad}a).
\]
The last inclusion clearly holds if we consider $\operatorname*{ad}a$ as an
operator on $\overline{I}$. Since $\operatorname*{ad}a$ has countable
spectrum, $\mathcal{E}_{0}(\operatorname*{ad}a)$ is closed. So one can
consider $\mathcal{E}_{0}(\operatorname*{ad}a)$ as a Banach space, whence
$\operatorname*{ad}a$ is quasinilpotent on $I$.

In any of the given cases, $I_{0}$ consists of Engel elements of $I$, so by
Theorem \ref{main1}, $I$ is $E$-solvable. Thus, if $I=$ $L^{\prime\prime}$,
$L$ is reducible by Corollary \ref{onetwo}, and if $I=L^{\prime}$, $L$ is
$E$-solvable by Lemma \ref{ideal}.
\end{proof}

Theorem \ref{main1} can be partially extended to the case of a Lie algebra
graded by an infinite group. We present below an appropriate version of such
an extension.

Let $\Gamma$ be an arbitrary (not necessarily finite) commutative group, and
let $L:=\sum_{\gamma\in\Gamma}L_{\gamma}$ be a $\Gamma$-subgraded Lie algebra.
Let $\Lambda$ be another commutative group and $M:=\sum_{\lambda\in\Lambda
}M_{\lambda}$ be a $\Lambda$-subgraded Lie algebra. We say that $M$ is a
\textit{subgraded} \textit{subalgebra} of $L$ if for every $\delta\in\Lambda$
there is $\gamma\in\Gamma$ such that
\[
M_{\delta}\subset L_{\gamma}.
\]
We say that $L$ is \textit{locally finitely subgraded} if every finite subset
of homogeneous elements of $L$ is contained in a $\Lambda$-subgraded
subalgebra $M$ for some finite group $\Lambda$.

Examples of locally finitely subgraded Lie algebras include those which are
graded with a locally finite group $\Gamma$.

\begin{theorem}
\label{locfin} Let $\Gamma$ be a not necessarily finite, commutative group,
and let a $\Gamma$-subgraded Lie algebra $L\subset\mathcal{K}(\mathcal{X})$ be
locally finitely subgraded. If every component $L_{\gamma}$ consist of Engel
(respectively, Volterra) elements of $L$ then $L$ is $E$-solvable
(respectively, is Volterra).
\end{theorem}

\begin{proof}
Let $a,b\in L$ be arbitrary. Then we may write $a=\sum_{i=1}^{n}a_{i}$, and
$b=\sum_{j=1}^{m}b_{j}$, where all $a_{i},b_{j}$ are homogeneous elements of
$L$. The set $\{a_{i}:1\leq i\leq n\}\cup\{b_{j}:1\leq j\leq m\}$ is contained
in a finitely subgraded Lie subalgebra $M$ of $L$, and it is clear that every
component of $M$ consists of Engel elements of $M$, whence $M$ is
triangularizable. In particular, $[a,b]$ is Volterra, whence by Theorem
\ref{redcri}, $L$ is $E$-solvable.

The case with Volterra components is proved similarly. Indeed, since, by
Corollary \ref{volt}, every finitely subgraded Lie algebra of compact
operators with Volterra components is Volterra, $a$ is a Volterra operator.
\end{proof}

\part{\label{__3related}$\mathbb{Z}_{2}$-Subgraded Lie Algebras, Lie Triple
Systems and Jordan Algebras}

The situation we treat in this part differs from one in Part \ref{__2graded}.
For a $\mathbb{Z}_{2}$-subgraded Lie algebra of compact operators
$L=L_{0}+L_{1}$, we impose a quasinilpotence condition to $L_{1}$ instead of
$L_{0}$. The main result is Theorem \ref{two}.

As an application, we establish that every Lie triple system or Jordan algebra
of Volterra operators is triangularizable, and also that every Jordan algebra
of compact operators containing a non-zero Volterra ideal is reducible.

\section{\label{_9z2}Lie Algebras Graded by $\mathbb{Z}_{2}$}

$\mathbb{Z}_{2}$-subgraded Lie algebras are of special interest. We have
already seen that if $L=L_{0}+L_{1}$ is a $\mathbb{Z}_{2}$-subgraded Lie
algebra of compact operators, and $L_{0}$ consists of Engel elements of $L$,
then $L$ is triangularizable. It turns out that if $L_{1}$ generates $L$ and
consists of Engel elements of $L$ then $L$ is triangularizable. As the example
\ref{e2} below shows, this does not extend to $\mathbb{Z}_{n}$-subgraded Lie
algebras for arbitrary $n$.

Note that if $L=L_{0}+L_{1}$ is a $\mathbb{Z}_{2}$-subgraded Lie algebra of
operators then
\[
L_{-1}=L_{1}%
\]
and so $L^{\prime\prime}$ is equal to $L_{1}+[L_{1},L_{1}]$. Recall that
$L^{\prime\prime}$ is an ideal of $L$.

\begin{lemma}
\label{findim2} Let $L=L_{0}+L_{1}$ be a $\mathbb{Z}_{2}$-subgraded Lie
algebra of operators on a finite dimensional space $\mathcal{X}$. If $L_{1}$
consists of Engel elements of $L$, then $L^{\prime\prime}$ is solvable. In
particular, if $L_{1}$ is non-scalar, then $L$ is reducible.
\end{lemma}

\begin{proof}
Let $I=L^{\prime\prime}$. Since $I$ is an ideal of $L$, it suffices to show
that $I$ is triangularizable. For this, it suffices to prove that $I$ is reducible.

Suppose, to the contrary, that $I$ is irreducible. Then setting $I_{0}%
=[L_{1},L_{1}]$ and $I_{1}=L_{1}$, we have that
\[
I=I_{0}+I_{1}%
\]
is a $\mathbb{Z}_{2}$-subgraded Lie algebra. Consider the following two cases.

\begin{case}
\label{caseone} $I_{1}$ consists of nilpotent operators.
\end{case}

Note that $I_{0}\cap I_{1}$ is an ideal of $I$, and that moreover, since it
consists of nilpotent operators, it is a Volterra ideal of $L$. If it is
non-zero then $I$ is reducible (by Theorem \ref{solid}), which would be a
contradiction. So we have that%
\[
I_{0}\cap I_{1}=0.
\]
In other words, $I$ is a $\mathbb{Z}_{2}$-graded Lie algebra.

Since $L_{1}$ consists of nilpotent operators, we have that
\[
\operatorname*{tr}(a^{2})=0
\]
for each $a\in L_{1}$. Since $\operatorname*{tr}(ab)=\frac{1}{2}%
\operatorname*{tr}((a+b)^{2}-a^{2}-b^{2})$, we see that
\[
\operatorname*{tr}(ab)=0
\]
for all $a,b\in L_{1}$. In other words, we have that
\[
\operatorname*{tr}(L_{1}L_{1})=0.
\]
Using this and $[L_{1},[L_{1},L_{1}]]\subset L_{1}$, we obtain that
\[
\operatorname*{tr}(I_{0}I_{0})=\operatorname*{tr}([L_{1},L_{1}][L_{1}%
,L_{1}])=\operatorname*{tr}(L_{1}[L_{1},[L_{1},L_{1}]])=0.
\]
By Lemma \ref{cartan}, the Lie algebra $I_{0}$ is solvable. So if $I_{0}$ is
not commutative, Lemma \ref{L0-triang} implies that $I$ is reducible, which
would be a contradiction. Therefore, $I_{0}$ is commutative. Moreover, by
Lemma \ref{prime}, $I_{0}$ is not scalar.

Set $N=[I_{0},I_{1}]$. If $N$ is scalar then $N=0$ (since the operators in $N$
have zero trace), and $I_{0}$ is in the center of $I$, whence $I$ is reducible
by Lomonosov's Theorem \cite{Lom}, which would be a contradiction. So $N$ is
not scalar.

Since $[I_{0},I_{0}]=0$, one has
\[
\operatorname*{tr}(I_{0}N)=\operatorname*{tr}(I_{0}[I_{0},I_{1}%
])=\operatorname*{tr}(I_{1}[I_{0},I_{0}])=0.
\]
Furthermore, since $N\subset I_{1}=L_{1}$, we have that
\[
\operatorname*{tr}(I_{1}N)=0.
\]
Therefore we obtain that
\[
\operatorname*{tr}(IN)=0,
\]
and hence by Lemma \ref{carfin}, $I$ is reducible, which is a contradiction.

This proves that $I$ is, in fact, triangularizable.

\begin{case}
$I_{1}$ consists of Engel elements of $L$.
\end{case}

In this case we use the adjoint representation $\operatorname*{ad}%
=\operatorname*{ad}_{L}$ of $L$ on $L$. Indeed, $(\operatorname*{ad}%
L)^{\prime\prime}=\operatorname*{ad}I$, the first component
$(\operatorname*{ad}L)_{1}$ of $\operatorname*{ad}L$ equals
$\operatorname*{ad}I_{1}$, and hence consists of nilpotent operators. Then by
Case \ref{caseone}, $\operatorname*{ad}I$ is triangularizable, and hence is
solvable. This means that $I$ is solvable, whence it is triangularizable.

In either case, if $I$ is not scalar, $L$ is reducible by Theorem \ref{solid}.
\end{proof}

\begin{lemma}
\label{finrank2} Let $L$ be a $\mathbb{Z}_{2}$-subgraded Lie algebra of finite
rank operators. If $L_{1}$ consists of Engel elements of $L$, then
$L^{\prime\prime}$ is $E$-solvable. In particular, if $L_{1}$ is non-scalar,
then $L$ is reducible.
\end{lemma}

\begin{proof}
The lemma can be proved in the same way as Lemma \ref{finrank}, but with the
use of Lemma \ref{findim2} instead of Lemma \ref{finsubgraded}.
\end{proof}

\begin{lemma}
\label{onefin2} Let $L$ be a $\mathbb{Z}_{2}$-subgraded Lie algebra of compact
operators. If $L_{1}$ consists of Engel elements of $L$, and contains a
non-scalar finite rank operator, then $L$ is reducible.
\end{lemma}

\begin{proof}
The lemma follows from Lemma \ref{finrank2} in the same way that Lemma
\ref{onefin} followed from Lemma \ref{finrank}.
\end{proof}

\begin{theorem}
\label{two}Let $L=L_{0}+L_{1}$ be a $\mathbb{Z}_{2}$-subgraded Lie algebra of
compact operators. If $L_{1}$ is non-scalar, and consists of Engel elements of
$L$, then $L$ is reducible.
\end{theorem}

\begin{proof}
Suppose to the contrary, that $L$ is irreducible.

Assume first that $L_{0}$ and $L_{1}$ are closed. By Lemma \ref{onefin2},
$L_{1}$ doesn't contain any non-scalar finite rank operators. By Theorem
\ref{main1}, $L_{0}$ doesn't consist of Engel elements of $L$.

If there is $a\in L_{0}$ such that $\operatorname*{ad}a$ is not quasinilpotent
on $L_{1}$, then there is a non-zero isolated point $\lambda\in\sigma
((\operatorname*{ad}a)|L_{1})$, and a non-zero $x\in L_{1}$ such that
\begin{equation}
x\in\mathcal{E}_{\lambda}(\operatorname*{ad}a) \label{formula}%
\end{equation}
Then $x$ is a nilpotent finite rank operator in $L_{1}$ by \cite[Theorem
3]{Woj}, which would contradict our assumption. Thus $(\operatorname*{ad}%
L_{0})|L_{1}$ consists of quasinilpotent operators.

We claim that there is $a\in L_{0}$ such that $\operatorname*{ad}a$ is not
quasinilpotent on $L_{0}$. Indeed, otherwise $\operatorname*{ad}L_{0}$
consists of operators that are locally quasinilpotent on $L_{0}\cup L_{1}$ and
hence are quasinilpotent on $L_{0}+L_{1}$ (see the proof of Corollary
\ref{nilfinrank}), which would be a contradiction.

Therefore, there exists a non-zero isolated point $\lambda\in\sigma
((\operatorname*{ad}a)|L_{0})$, and a non-zero $x\in L_{0}$ such that
(\ref{formula}) holds. This implies that $x$ is a nilpotent finite rank
operator in $L_{0}$. Since $L_{1}$ does not contain non-scalar finite rank
operators, we have that
\[
\lbrack x,L_{1}]=0.
\]
Moreover, we must also have that
\[
\lbrack x,L^{\prime\prime}]=0.
\]
Then, since $x$ commutes with $L^{\prime\prime}$, $xL^{\prime\prime}$ consists
of nilpotent finite rank operators, and we obtain that
\[
\operatorname*{tr}(xL^{\prime\prime})=0.
\]
Let $I=\{b\in L\cap\mathcal{F}(\mathcal{X}):\operatorname*{tr}(bL^{\prime
\prime})=0\}$. Since $x\in I$, the set $I$ is non-scalar. Since $[L,L^{\prime
\prime}]\subset L^{\prime\prime}$, it follows from
\[
\operatorname*{tr}([I,L]L^{\prime\prime})=\operatorname*{tr}(I[L,L^{\prime
\prime}])=0
\]
that $I$ is an ideal of $L$. Let $N=I\cap L^{\prime\prime}$. Then $N$ is an
ideal of $L$, and
\[
\operatorname*{tr}(NN)=0.
\]
By Theorem \ref{carMatt}, $[N,N]$ consists of nilpotent operators, whence $N$
is $E$-solvable by Theorem \ref{redcri}. If $N$ is non-scalar then $L$ is
reducible by Theorem \ref{solid}, which would be a contradiction.

So $N$ is scalar. Since $I$ and $L^{\prime\prime}$ are non-scalar ideals of
$L$ with scalar intersection, $L$ is reducible by Theorem \ref{twoideals}, a contradiction.

This proves that $L$ is reducible.

In the general case, consider the $\mathbb{Z}_{2}$-subgraded Lie algebra $M$
with $M_{0}=\overline{L_{0}}$ and $M_{1}=\overline{L_{1}}$. It is easy to
check that $M_{1}$ consists of Engel elements of $M$ and that $M^{\prime}$ is
non-scalar. Then the previous argument applied to $M$ shows that $M$ is
reducible, whence $L$ is reducible.
\end{proof}

We collect our results on $\mathbb{Z}_{2}$-subgraded Lie algebras in the
following theorem.

\begin{theorem}
\label{2theor} Let $L$ be a $\mathbb{Z}_{2}$-subgraded Lie algebra of compact
operators. Then

\begin{enumerate}
\item[(i)] $L$ is reducible if $L_{1}$ is non-scalar, and if one of the
following conditions holds.

\begin{itemize}
\item $L_{1}$ consists of Engel elements of $L^{\prime\prime}$.

\item $[L_{1},L_{1}]$ consists of Engel elements of $L^{\prime\prime}$.
\end{itemize}

\item[(ii)] $L$ is $E$-solvable if one of the following conditions holds.

\begin{itemize}
\item $L_{0}$ consists of Engel elements of $L$.

\item $L_{1}$ generates $L$ (as a Lie algebra) and consists of Engel elements
of $L$.

\item $[L_{0},L_{0}]$ and $[L_{1},L_{1}]$ consist of Engel elements of $L$.
\end{itemize}
\end{enumerate}
\end{theorem}

\begin{proof}
In both first cases of (ii) it suffices to show that $L$ is reducible. These
statements follow from Theorems \ref{main1} and \ref{two}. The last case of
(ii) follows from Theorem \ref{lieset}.

In the case (i), applying results of (ii) to $L^{\prime\prime}$, we obtain
that $L^{\prime\prime}$ is $E$-solvable. Since $L^{\prime\prime}$ is a
non-scalar ideal of $L$, $L$ is reducible by Theorem \ref{solid}.
\end{proof}

\begin{remark}
$\emph{Recall}$\emph{ }$\emph{that}$ each Volterra operator is an Engel
element of every operator Lie algebra which contains it. $\emph{Thus,}$
$\emph{it}$ $\emph{follows}$ $\emph{from}$ $\emph{Theorem}$ $\ref{2theor}$
$\emph{that}$

\begin{enumerate}
\item[(i)] If $L_{1}$ is Volterra and non-zero then $L$ is reducible.

\item[(ii)] If $L_{1}$ is Volterra and generates $L$ (as a Lie algebra) then
$L$ is triangularizable.
\end{enumerate}
\end{remark}

\section{\label{10triple}Operator Lie Triple Systems}

A\ subspace $M$ of $\mathcal{B}(\mathcal{X})$ is a \textit{Lie triple system}
if it is closed under the Lie triple product $[a,[b,c]]$ for all $a,b,c\in M$.
It is clear that every Lie algebra of operators\ is a Lie triple system.

The following lemma shows that there is a natural embedding of Lie triple
systems into $\mathbb{Z}_{2}$-subgraded Lie algebras.

\begin{lemma}
\label{Jor-grad} Let $M$ be a Lie triple system in $\mathcal{B}(\mathcal{X})$.
Set $L_{0}=[M,M]$ and $L_{1}=M$. Then $L=L_{0}+L_{1}$ is a $\mathbb{Z}_{2}%
$-subgraded Lie algebra.
\end{lemma}

\begin{proof}
The inclusion $[L_{1},L_{1}]\subset L_{0}$ is evident; in fact,
\[
\lbrack L_{1},L_{1}]=L_{0}.
\]
The inclusion
\[
\lbrack L_{0},L_{1}]\subset L_{1}%
\]
is also clear; it follows from the definition of a Lie triple system. Finally,
we have
\begin{align*}
\lbrack L_{0},L_{0}]  &  =[L_{0},[L_{1},L_{1}]]\subset\lbrack\lbrack
L_{0},L_{1}],L_{1}]+[[L_{1},[L_{0},L_{1}]]\\
&  \subset\lbrack L_{1},L_{1}]=L_{0}\text{.}%
\end{align*}

\end{proof}

Given a subset $M$ of $\mathcal{B}(\mathcal{X})$, we denote by $\mathcal{L}%
(M)$ the Lie algebra generated by $M$. If $M$ is a Lie triple system then it
is clear that%
\[
\mathcal{L}(M)=M+[M,M].
\]

\begin{theorem}
\label{tripvolt} A Lie triple system $M$ of Volterra operators is triangularizable.
\end{theorem}

\begin{proof}
Let $L=\mathcal{L}(M)$. It suffices to prove that the $\mathbb{Z}_{2}%
$-subgraded Lie algebra $L$ is triangularizable. Since $L_{1}=M$ and $M$ is
Volterra, this follows by Theorem \ref{2theor}.
\end{proof}

For a subset $M$ of a Lie algebra $L$, put
\[
M^{\left[  1\right]  }=M\quad\text{and}\quad M^{\left[  k+1\right]  }=\left\{
[a,b]:a\in M,b\in M^{\left[  k\right]  }\right\}  .
\]
It is clear that $\cup M^{\left[  k\right]  }$ is a Lie multiplicative set in
$L$. For a subspace $M$ of $L$ and $n>1$, we say that $M$ is a \textit{Lie}
$n$-\textit{product system} in $L$ if
\[
M^{\left[  n\right]  }\subset M.
\]

The following example shows that the result of Theorem \ref{tripvolt} does not
extend to Lie $n$-product systems of Volterra operators for $n=5$.

\begin{example}
\label{e2} Let $a=\left(
\begin{array}
[c]{ccc}%
0 & 1 & 0\\
0 & 0 & -1\\
0 & 0 & 0
\end{array}
\right)  $, $b=\left(
\begin{array}
[c]{ccc}%
0 & 0 & 0\\
1 & 0 & 0\\
0 & 1 & 0
\end{array}
\right)  $, and let $M$ be the linear space generated by $a$ and $b$. \emph{It
is not difficult to check that} $M$ is an irreducible Lie $5$-product system
of nilpotent finite rank operators.

\emph{As a consequence,} there is an irreducible $\mathbb{Z}_{4}$-subgraded
Lie algebra $L$ of finite rank operators such that the component $L_{1}$
consists of nilpotents and generates $L$ as a Lie algebra. \emph{Indeed, let
}$L=\mathcal{L}(M)$. \emph{Then}
\[
L=M^{\left[  1\right]  }+\operatorname*{span}M^{\left[  2\right]
}+\operatorname*{span}M^{\left[  3\right]  }+\operatorname*{span}M^{\left[
4\right]  }%
\]
\emph{in virtue of} $M^{\left[  5\right]  }\subset M$. \emph{Take
}$\operatorname*{span}M^{\left[  4\right]  }$ \emph{as the zero component of}
$L$ \emph{and} $\operatorname*{span}M^{\left[  i\right]  }$ \emph{as the}
$i$-\emph{th} \emph{component of} $L$ \emph{for} $i=1,2,3$. \emph{It is easy
to check that} $L$ \emph{is} $\mathbb{Z}_{4}$-\emph{subgraded}.
\end{example}

\section{\label{11jordan}Operator Jordan Algebras}

Recall that a subspace $J$ of $\mathcal{B}(\mathcal{X})$ is a \textit{Jordan
algebra} if it is closed under the Jordan product
\[
a\circ b=ab+ba
\]
for all $a,b\in J$. The equality%
\begin{equation}
\lbrack a,[b,c]]=(a\circ b)\circ c-(a\circ c)\circ b \label{joriden}%
\end{equation}
shows that \textit{every Jordan algebra is also a Lie triple system}.

\begin{corollary}
\label{jordvolt} A Jordan algebra of Volterra operators is triangularizable.
\end{corollary}

\begin{proof}
Follows from Theorem \ref{tripvolt}.
\end{proof}

For the case of a Jordan algebra of Shatten operators on a Hilbert space this
result was obtained in \cite{Ken}. It was a starting point of the present work.

Let $J$ be a Jordan algebra. Recall that a \textit{Jordan ideal }$I$ of $J$ is
a subspace of $J$ such that%
\[
J\circ I\subset I.
\]
Here, for Lie algebras $L$ and $N$, the designation $N\vartriangleleft L$
means that $N$ is an ideal of $L$.

\begin{lemma}
\label{jorideals}Let $J$ be a Jordan algebra in $\mathcal{B}(\mathcal{X})$ and
let $I$ be a Jordan ideal of $J$. Define
\[
\mathcal{L}(J,I)=I+[J,I].
\]
Then $\mathcal{L}(J,I)$ is a Lie algebra, and%
\[
\mathcal{L}(I)\vartriangleleft\mathcal{L}(J,I)\vartriangleleft\mathcal{L}(J)
\]
is a series of Lie ideals, where $\mathcal{L}(I)$ and $\mathcal{L}(J)$ are Lie
algebras generated by $I$ and $J$, respectively.
\end{lemma}

\begin{proof}
Recall that $\mathcal{L}(J)=J+[J,J]$ and $\mathcal{L}(I)=I+[I,I]$. Since $I$
is an ideal of $J$, the inclusions
\[
\lbrack\lbrack J,I],J]\subset I\text{ and }[[J,J],I]\subset I
\]
easily follow from (\ref{joriden}). Then, applying the Jacobi identity, we
obtain that%
\[
\lbrack\lbrack J,J],[J,I]]\subset\lbrack\overset{~\quad\subset J}%
{\overbrace{[[J,J],J]}},I]+[J,\overset{~\quad\subset I}{\overbrace{[[J,J],I]}%
}]\subset\lbrack J,I],
\]
and%
\[
\text{ }[[J,I],[I,I]]\subset\lbrack\underset{~\quad\subset I}{\underbrace
{[[J,I],I]}},I]+[I,\underset{~\quad\subset I}{\underbrace{[[J,I],I]}}%
]\subset\lbrack I,I].
\]
It is now easy to verify that%
\[
\lbrack\mathcal{L}(J,I),\mathcal{L}(I)]\subset\mathcal{L}(I)\text{ and
}[\mathcal{L}(J),\mathcal{L}(J,I)]\subset\mathcal{L}(J,I)\text{.}%
\]

\end{proof}

\begin{theorem}
\label{jordan} A Jordan algebra $J$ of compact operators with a non-zero
Volterra ideal $I$ is reducible.
\end{theorem}

\begin{proof}
As in Lemma \ref{jorideals}, we have the series of ideals of Lie algebras
\[
\mathcal{L}(I)\vartriangleleft\mathcal{L}(J,I)\vartriangleleft\mathcal{L}(J).
\]
By Corollary \ref{jordvolt}, $I$ is triangularizable. Then $\mathcal{L}(I)$ is
triangularizable and therefore is $E$-solvable by Theorem \ref{redcri}. Let
$R$ be the largest $E$-solvable ideal of $\mathcal{L}(J,I)$. Since
$\mathcal{L}(I)$ is an $E$-solvable ideal of $\mathcal{L}(J,I)$, we obtain
that
\[
I\subset\mathcal{L}(I)\subset R.
\]
In particular, $R$ contains non-scalar elements. By Theorem \ref{inner}, $R$
is inner-characteristic. This means that
\[
\lbrack a,R]\subset R
\]
for every $a\in\mathcal{B}(\mathcal{X})$ with $[a,\mathcal{L}(J,I)]\subset
\mathcal{L}(J,I)$. Since
\[
\lbrack\mathcal{L}(J),\mathcal{L}(J,I)]\subset\mathcal{L}(J,I),
\]
this implies that
\[
\lbrack\mathcal{L}(J),R]\subset R.
\]
In other words, $R$ is an ideal of Lie algebra $\mathcal{L}(J)$. Since, by
Theorem \ref{solid}, every Lie algebra of compact operators with a non-scalar
$E$-solvable ideal is reducible, we obtain that $\mathcal{L}(J)$ is reducible,
and hence that $J$ is reducible.
\end{proof}

\end{document}